\documentclass{amsart}
\usepackage{oldlfont,amsfonts,amssymb}
\newtheorem{theorem}{Theorem}[section]
\newtheorem{lemma}[theorem]{Lemma}

\theoremstyle{definition}
\newtheorem{definition}[theorem]{Definition}
\newtheorem{example}[theorem]{Example}

\newtheorem{problem}[theorem]{Problem}
\theoremstyle{remark}
\newtheorem{remark}[theorem]{Remark}

\numberwithin{equation}{section}

\begin{document}

\title{On Bures distance over Standard Form $vN$-Algebras}

\author{Peter M. Alberti}
\address{Institute of Theoretical Physics\\
University of Leipzig\\
Augustusplatz 10, D-04109 Leipzig, Germany}
\email{Peter.Alberti@itp.uni-leipzig.de}
\thanks{The paper is the completed version of some parts
of the lectures on Bures geometry which were held by the first
author at `Graduiertenkolleg Quantenfeldtheorie', in spring term
2000.}

\author{Gregor Peltri}
\address{Naturwissenschaftlich-Theoretisches Zentrum\\University of Leipzig\\
Augustusplatz 10, D-04109 Leipzig, Germany}
\email{Gregor.Peltri@itp.uni-leipzig.de}
\thanks{The second author was supported by a `Doktorandenf\"orderplatz' at ITP}

\subjclass{46L89, 46L10, 58B20}
\date{}

\commby{}

\dedicatory{}

\begin{abstract}
In case of standard form $vN$-algebras, the Bures distance is the natural distance
between the fibres of
implementing vectors at normal positive linear forms.
Thereby, it is well-known that to each two normal positive linear
forms implementing vectors exist such that the Bures distance
is attained by the metric distance of the implementing vectors in question.
We discuss to which extent this can remain true if a vector in one of the fibres is considered
as fixed. For each nonfinite
algebra, classes of counterexamples are given and situations are analyzed where the latter type of
result must fail. In the course of the paper, an account of those facts and notions is given,  
which can be taken as a useful minimum of basic ${\mathsf C}^*$-algebraic tools needed
in order to efficiently develop the fundamentals of Bures geometry over standard form $vN$-algebras.
\end{abstract}
\maketitle
\section{Basic settings and results}\label{bas}
\subsection{Definitions and conventions}\label{bas1}
Throughout the paper, a distance function $d_{\mathrm{B}}$ on the
positive cone $M_+^*$ of the bounded linear forms $M^*$ over a unital
${\mathsf C}^*$-algebra $M$ will be considered. For normal states on
a ${\mathsf W}^*$-algebra $d_{\mathrm{B}}$ agrees with
the Bures distance function \cite{Bure:69}. Therefore,
henceforth also $d_{\mathrm{B}}$ will be referred to as Bures distance function.
We start by defining $d_{\mathrm{B}}(M|\nu,\varrho)$
between $\nu,\varrho\in M_+^*$.
\begin{definition}\label{budi}
$\ \ d_{\mathrm{B}}(M|\nu,\varrho)=\inf_{\{\pi,{\mathcal K}\},\varphi\in
{\mathcal S}_{\pi,M}(\nu),\psi\in {\mathcal S}_{\pi,M}(\varrho)}
\|\psi-\varphi\|\,.$
\end{definition}
Instead of $d_{\mathrm{B}}(M|\nu,\varrho)$ the notation
$d_{\mathrm{B}}(\nu,\varrho)$ will often be used.
For unital $^*$-representation $\{\pi,{\mathcal K}\}$ of $M$ on a Hilbert space $\{{\mathcal K},
\langle\cdot,\cdot\rangle\}$
and for $\mu\in M_+^*$ we let
\begin{equation}\label{faser}
{\mathcal S}_{\pi,M}(\mu)=\{\chi\in {\mathcal K}:\mu(\cdot)=
\langle\pi(\cdot)\chi,\chi\rangle\}.
\end{equation}
In case of ${\mathcal S}_{\pi,M}(\mu)\not=\emptyset$ this set will be referred
to as $\pi${\em{-fibre of $\mu$}}. The above infimum extends
over all $\pi$ relative to which both $\pi$-fibres exist and,
within each
such representation, $\varphi$ and $\psi$ may be varied
through all of ${\mathcal S}_{\pi,M}(\nu)$ and ${\mathcal S}_{\pi,M}(\varrho)$, respectively.
The scalar product ${\mathcal K}\times {\mathcal K}
\ni\{\chi,\eta\}\,\longmapsto\,\langle\chi,\eta\rangle\in
{\mathbb{C}}$ on the representation Hilbert space by convention is supposed
to be linear with respect to the first argument $\chi$, and antilinear
in the second argument $\eta$, and maps into the complex
field ${\mathbb{C}}$. Let ${\mathbb{C}}\ni z\mapsto \bar{z}$ be the complex conjugation,
and be $\Re z$ and $|z|$ the real part and absolute value of $z$, respectively. The norm of $\chi\in
{\mathcal K}$ is given by $\|\chi\|=\sqrt{\langle\chi,\chi\rangle}$.
For the relating operator and ${\mathsf C}^*$-algebra theory, the reader is referred to the standard
monographs,
e.g.~\cite{Dixm:64,Saka:71,KaRi:83}.

For both the ${\mathsf C}^*$-norm of an element
$x\in M$ as well as for the operator norm of a concrete bounded linear
operator $x\in {\mathsf B}({\mathcal K})$ the same notation
$\|x\|$ will be used, and the involution ($^*$-operation)
respectively the taking of the hermitian conjugate of an element $x$ is
indicated by the transition $x\,\longmapsto\,x^*$. The notions of hermiticity and positivity
for elements are defined as usual in ${\mathsf C}^*$-algebra theory, and
$M_{\mathrm{h}}$ and $M_+$ are the hermitian and positive elements of $M$,
respectively. The null
and the unit element/operator in $M$ and ${\mathsf B}({\mathcal
K})$ will be denoted by ${\mathbf 0}$ and ${\mathbf 1}$.
For notational purposes mainly, in short recall some
fundamentals relating (bounded) linear forms which subsequently
might be of concern in context of Definition \ref{budi}. Recall that the
topological dual space $M^*$ of $M$ is the
set of all those linear functionals (linear forms) which are continuous with respect to the
operator norm topology. Equipped with the dual norm
$\|\cdot\|_1$, which is
given by $\|f\|_1=\sup\{|f(x)|\,:\,x\in M,\,\|x\|\leq 1\}$
and which is referred to as the functional norm, $M^*$ is a Banach space.
For each given $f\in M^*$, the hermitian
conjugate functional $f^*\in {M}^*$ is defined by $f^*(x)=\overline{f(x^*)}$, for
each $x\in M$. Remind that $f\in M^*$ is hermitian
if $f=f^*$ holds, and $f$ is termed positive if $f(x)\geq 0$ holds,
for each $x\in M_+$. Relating positivity, the basic fact is that
a bounded linear form over $M$ is positive
if, and only if, $\|f\|_1=f({\mathbf 1})$.

For each $f\in M_+^*$, there exists a cyclic $^*$-representation
$\pi_f$ of $M$ on some Hilbert space ${\mathcal K}_f$, with cyclic vector
$\varOmega\in {\mathcal K}_f$, and obeying $f(x)=\langle \pi_f(x)\varOmega,
\varOmega\rangle$, for all $x\in M$ (Gelfand-Neumark-Segal theorem).
Considering that construction in the special case
with $f=\nu+\varrho$ will provide a unital $^*$-representation $\pi=\pi_f$
where the $\pi$-fibres of $\nu$ and $\varrho$ both exist (we omit the details, all of which are standard). Thus Definition
\ref{budi} makes sense, in any case of $\nu,\varrho\in M_+^*$.

In conjunction with the Bures distance $d_{\mathrm B}$ there appears the functor $P$ of the
($^*$-algebraic) {\em transition probability} \cite{Uhlm:76}.
For given ${\mathsf C}^*$-algebra $M$ and positive
linear forms $\nu,\varrho\in M_+^*$ the definition reads as follows\,:
\begin{definition}\label{genprob.2}
\ \ $P_M(\nu,\varrho)=
\sup_{\{\pi,{\mathcal K}\},\varphi\in
{\mathcal S}_{\pi,M}(\nu),\psi\in {\mathcal S}_{\pi,M}(\varrho)}
|\langle \psi,\varphi\rangle|^2\,.$
\end{definition}
The range of variables in the supremum is the same as in Definition \ref{budi}.
With the help of $P_M$ the following formula for $d_{\mathrm{B}}$ is obtained\,:
\begin{equation}\label{pcont.1}
d_{\mathrm{B}}(M|\nu,\varrho)^2=\biggl\{\|\nu\|_1-
\sqrt{P_M(\nu,\varrho)}\biggr\} +
\biggl\{\|\varrho\|_1-\sqrt{P_M(\nu,\varrho)}\biggr\}\,.
\end{equation}
Thus, properties of $d_{{\mathsf B}}$ can be obtained from
properties of $P$, and vice versa.
\begin{remark}\label{orig}
In \cite{Bure:69}, in generalizing from the commutative case which had been studied extensively
in \cite{Kaku:48}, expressions $d$ and $\rho$ instead of $d_{\mathrm{B}}$ and $\sqrt{P}$
were considered for normal states
on a ${\mathsf W}^*$-algebra, but with the infima
extending over all faithful representations $\{\pi,{\mathcal K}\}$ where both fibres exist, accordingly.
\end{remark}
\subsection{Some basic results on Bures distance}\label{bas2}
Let $\nu,\varrho\in M_+^*$, and be $\{\pi,{\mathcal K}\}$ a unital $^*$-representation of $M$
such that the $\pi$-fibres of $\nu$ and $\varrho$ both exist. Suppose $\varphi\in {\mathcal S}_{\pi,M}(\nu)$ and
$\psi\in {\mathcal S}_{\pi,M}(\varrho)$. Let
$$\pi(M)^{\,\prime}=\{z\in {\mathsf B}({\mathcal K}): z\pi(x)=\pi(x)z,\,\forall x\in M\}$$
be the commutant $vN$-algebra of $\pi(M)$, and be ${\mathcal U}(\pi(M)^{\,\prime})$ the group
of unitary operators of $\pi(M)^{\,\prime}$.
Define a linear form $h_{\psi,\varphi}^\pi$ as follows\,:
\begin{equation}\label{hform}
\forall \,z\in\pi(M)^{\,\prime}\,:\  h_{\psi,\varphi}^\pi(z)=\langle z\psi,\varphi\rangle\,.
\end{equation}
For $\chi\in {\mathcal K}$ let orthoprojections
$p_\pi(\chi),\,p_\pi^{\prime}(\chi)$ be defined as the orthoprojections projecting from ${\mathcal K}$
onto the the closed
linear subspaces $\overline{\pi(M)^{\,\prime}\chi}$ and
$\overline{\pi(M)\chi}$, respectively.
It is standard that $p_\pi(\chi)\in \pi(M)^{\,\prime\prime}$ and
$p_\pi^{\,\prime}(\chi)\in \pi(M)^{\,\prime}$ hold, with the double commutant
$vN$-algebra $$\pi(M)^{\,\prime\prime}=
\bigl(\pi(M)^{\,\prime}\bigr)^{\,\prime}$$ of $\pi(M)$. By the Kaplansky--von\,Neumann theorem
$\pi(M)$ is strongly dense within $\pi(M)^{\,\prime\prime}$, and therefore
one always has
\begin{equation}\label{bas3a}
\overline{\pi(M)\chi}=\overline{\pi(M)^{\,\prime\prime}\chi}\,,
\end{equation}
which is useful to know. Relating \eqref{faser}, for each
$\chi\in {\mathcal S}_{\pi,M}(\mu)$ the following is true\,:
\begin{equation}\label{bas4}
{\mathcal S}_{\pi,M}(\mu)=\bigl\{v\chi:\,v^*v=p_\pi^{\prime}(\chi),\,v\in \pi(M)^{\,\prime}\bigr\}\,.
\end{equation}
Recall some useful facts about \eqref{hform}.
Assume $\varphi\in {\mathcal S}_M(\nu)$,
$\psi\in {\mathcal S}_M(\varrho)$, and be
\begin{subequations}\label{sub}
\begin{equation}\label{sub0}
h_{\psi,\varphi}^\pi=
\bigl|h_{\psi,\varphi}^\pi\bigr|\bigl((\cdot)v_{\psi,\varphi
}^\pi\bigr)
\end{equation}
the polar decomposition of the normal linear form $h_{\psi,\varphi}^\pi$.
According to the polar decomposition
theorem for normal linear forms in $vN$-algebras, a partial isometry
$v=v_{\psi,\varphi}^\pi$ in and a normal positive linear form
$g=\bigl|h_{\psi,\varphi}^\pi\bigr|$ over $\pi(M)^{\,\prime}$ and obeying
$h_{\psi,\varphi}^\pi=g((\cdot)v)$ exist, and both
are unique subject to $v^*v=s(g)$, with the support orthoprojection  $s(g)$
of $g$. Thus especially
\begin{subequations}\label{subsub1}
\begin{equation}\label{sub0.5}
\bigl|h_{\psi,\varphi}^\pi\bigr|=h_{\psi,\varphi}^\pi\bigl((\cdot)v_{\psi,\varphi}^{\pi*}\bigr)=
\bigl\langle(\cdot)v_{\psi,\varphi}^{\pi*}\psi,\varphi\bigr
\rangle\,,
\end{equation}
from which in view of the above the following can be obtained:
\begin{equation}\label{sub1}
v_{\psi,\varphi}^{\pi*}v_{\psi,\varphi}^\pi=
s\bigl(\bigl|h_{\psi,\varphi}^\pi\bigr|\bigr)\leq p_\pi^\prime(\varphi)\,.
\end{equation}
\end{subequations}
Analogously, by polar decomposition of $h_{\psi,\varphi}^{\pi*}=h_{\varphi,\psi}^\pi$,
one has $v_{\varphi,\psi}^\pi=v_{\psi,\varphi}^{\pi*}$ and
\begin{subequations}\label{subsub2}
\begin{equation}\label{sub1.5}
 \bigl|h_{\psi,\varphi}^{\pi*}\bigr|=h_{\varphi,\psi}^\pi\bigl((\cdot)v_{\psi,\varphi}^{\pi}\bigr)=
\bigl\langle(\cdot)v_{\psi,\varphi}^{\pi}\varphi,\psi\bigr
\rangle\,,
\end{equation}
from which analogously
\begin{equation}\label{sub2}
v_{\psi,\varphi}^\pi v_{\psi,\varphi}^{\pi*}=
s\bigl(\bigl|h_{\psi,\varphi}^{\pi*}\bigr|\bigr)\leq p_\pi^\prime(\psi)
\end{equation}
\end{subequations}
is obtained. In particular, with the Murray-von Neumann equivalence `$\sim$' one has
\begin{equation}\label{sub3}
s\bigl(\bigl|h_{\psi,\varphi}^\pi\bigr|\bigr)\sim
s\bigl(\bigl|h_{\psi,\varphi}^{\pi*}\bigr|\bigr)\,.
\end{equation}
\end{subequations}
We are going to comment now on some of the fundamentals of Bures geometry.
\subsubsection{Basic algebraic facts on Bures distance}\label{basfa}
From \cite[{\sc{Corollary 1}},\,{\sc{Corollary 2}},\,eqs.\,(5),\,(6) and {\sc{Theorem 3}}]{Albe:83}
the following facts are known\,:
\begin{lemma}\label{bas3}
$$\sqrt{P_M(\nu,\varrho)}=\bigl\| h_{\psi,\varphi}^\pi \bigr\|_1=
\sup_{u\in {\mathcal U}(\pi(M)^{\,\prime})}
\bigl| h_{\psi,\varphi}^\pi(u)\bigr|=\inf_{x\in M_+,\,\text{\rm{invertible}}} \sqrt{\nu(x)\varrho(x^{-1})}\,.$$
\end{lemma}
Each unitary orbit ${\mathcal U}(\pi(M)^{\,\prime})\chi$ is a
special subset of the $\pi$-fibre \eqref{bas4}. Thus,
the first two items of the following result at once can be seen from Lemma \ref{bas3}.
\begin{theorem}\label{bas5}
Suppose $\{\pi,{\mathcal K}\}$ is a unital $^*$-representation such that the $\pi$-fibres of
$\nu,\varrho\in M_+^*$ exist. For given $\varphi\in {\mathcal S}_{\pi,M}(\nu)$ and
$\psi\in {\mathcal S}_{\pi,M}(\varrho)$ the following hold\textup{:}
\begin{enumerate}
\item\label{bas51}
$d_{\mathrm{B}}(M|\nu,\varrho)=\inf_{u\in {\mathcal U}(\pi(M)^{\,\prime})}
\|u\psi-\varphi\|$\,;
\item\label{bas52}
$d_{\mathrm{B}}(M|\nu,\varrho)=\inf_{\psi'\in {\mathcal S}_{\pi,M}(\varrho)}
\|\psi'-\varphi\|$\,;
\item\label{bas53}
$d_{\mathrm{B}}(M|\nu,\varrho)=\|\psi-\varphi\|\ \Longleftrightarrow\ h_{\psi,\varphi}^\pi\geq 0$.
\end{enumerate}
\end{theorem}
Relating \eqref{bas53}, note that $h_{\psi,\varphi}^\pi\geq 0$ is equivalent
to $\langle\psi,\varphi\rangle=h_{\psi,\varphi}^\pi({\mathbf 1})=\bigl\| h_{\psi,\varphi}^\pi \bigr\|_1$.
The latter is equivalent to  $d_{\mathrm{B}}(M|\nu,\varrho)=\|\psi-\varphi\|$, by \eqref{pcont.1} and
Lemma \ref{bas3}.
\begin{remark}\label{rem1}
\begin{enumerate}
\item\label{rem10}
By Theorem \ref{bas5}, $d_{\mathrm{B}}$ can be calculated in each
$\pi$ where both fibres exist. Thus, on a ${\mathsf W}^*$-algebra $M$ and
for normal states,
the universal ${\mathsf W}^*$-representation can be used.
Hence, $d_{{\mathrm{B}}}$ and $\sqrt{P}$
argee with $d$ and $\rho$ of \cite{Bure:69} there, see Remark \ref{orig}.
\item\label{rem11}
For $\mu\in M_+^*$ and $\pi$ with ${\mathcal S}_{\pi,M}(\mu)\not=\emptyset$,
there is a unique normal positive $\mu_\pi$ on
$\pi(M)^{\,\prime\prime}$, with $\mu=\mu_\pi\circ\pi$.
By \eqref{bas3a} and \eqref{bas4} then
${\mathcal S}_{\mathrm{id},\pi(M)^{\,\prime\prime}}(\mu_\pi)={\mathcal S}_{\pi,M}(\mu)$,
with the trivial representation $\mathrm{id}$ of $\pi(M)^{\,\prime\prime}$ on ${\mathcal K}$.
\item\label{rem12}
For a $vN$-algebra $M$ and $\pi=\mathrm{id}$, the
notation of
the subscript/superscript `$\pi$' will
be omitted, and then ${\mathcal S}_M(\nu)$, $h_{\psi,\varphi}$, $v_{\psi,\varphi}$ and $p(\varphi)$, $p^\prime(\varphi)$
instead of ${\mathcal S}_{\pi,M}(\nu)$, $h_{\psi,\varphi}^\pi$, $v_{\psi,\varphi}^\pi$ and
$p_\pi(\varphi)$, $p_\pi^\prime(\varphi)$, with $\pi=\mathrm{id}$, respectively, will be in use and
${\mathcal S}_M(\nu)\not=\emptyset$ will be referred to as fibre of $\nu$, simply.
\end{enumerate}
\end{remark}
There are examples where Lemma \ref{bas3} can be made more explicit. Let $\tau$ be
a (lower semicontinuous)
trace on $M_+$ (see e.g.~in \cite[6.1.]{Dixm:64}). Then one has $^*$-ideals
${\mathcal L}^2(M,\tau)=\{x\in M:\tau(x^*x)<\infty\}$ and
${\mathcal L}^1(M,\tau)=\{x\in M:\,x=y^*z,\ y,z\in {\mathcal L}^2(M,\tau)\}$.
It is known that ${\mathcal L}^1(M,\tau)$
is the complex linear span of the hereditary positive cone
${\mathcal L}^1(M,\tau)_+=\{x\in M_+:\,\tau(x)<\infty\}$, and thus $\tau$ uniquely
extends to an {\em invariant} positive linear form $\tilde{\tau}$ onto
${\mathcal L}^1(M,\tau)$. That is, $\tilde{\tau}(xy)=\tilde{\tau}(yx)$,
for either $x,y\in {\mathcal L}^2(M,\tau)$, or $x\in M$ and
$y\in {\mathcal L}^1(M,\tau)$. By uniqueness, $\tau$ can be identified
with $\tilde{\tau}$, and then
is referred to as {\em trace $\tau$ on $M$}. In this sense, for $z\in {\mathcal L}^2(M,\tau)$ and
$y\in {\mathcal L}^1(M,\tau)$, we have a
unique positive linear form $\tau^z$ and bounded linear form $\tau_y$ on $M$, with
$\tau^z(x)=\tau(z^* xz)$ and $\tau_y(x)=\tau(yx)$, respectively, for each $x\in M$.

Let $I(M,\tau)=\{x\in M:\,\tau(x^*x)=0\}$. This is a $^*$-ideal
in ${\mathcal L}^2(M,\tau)$. Consider the
completion $L^2(M,\tau)$ of the space of
equivalence classes $\eta_x$ modulo $I(M,\tau)$ under the inner product
$\langle \eta_x, \eta_y\rangle=\tau(y^* x)$, $x,y\in {\mathcal L}^2(M,\tau)$. By standard arguments
the latter is well-defined and extends to a scalar product
on $L^2(M,\tau)$, which then is a Hilbert space. Let $\pi: M\in x\,\longmapsto\,\pi(x)$ and
$\pi^{\prime}(x): M\in x\,\longmapsto\,\pi^{\prime}(x)$ be the
$^*$-representations of $M$ over $L^2(M,\tau)$
which can be uniquely given through $\pi(x)\eta_y=\eta_{xy}$ and
$\pi^{\prime}(x)\eta_y=\eta_{yx}$, respectively, for each
$y\in {\mathcal L}^2(M,\tau)$ and all $x\in M$.
Then, for each $z\in {\mathcal L}^2(M,\tau)$ and $x\in M$ one has
$\tau^z(x)=\langle \pi(x)\eta_{z},\eta_{z}\rangle$. Note that $\pi^{\prime}(M)\subset
\pi(M)^{\,\prime}$. In carefully analyzing $\pi$ and $\pi'$ in case of a $vN$-algebra $M$, and
in applying Lemma \ref{bas3} with the
mentioned $^*$-representation $\pi$ and using that a normal trace is
lower semicontinuous, one can prove the following important special cases.
\begin{example}\label{ex1}
Let $M$ be a $vN$-algebra, with normal trace $\tau$. For $x,y\in {\mathcal L}^2(M,\tau)$ and
$a=|x^*|^2=xx^*$, $c=|y^*|^2=yy^*$, one has $\tau^x=\tau_a$ and $\tau^y=\tau_b$, and 
the following formulae hold\,:
\begin{subequations}\label{bei}
\begin{equation}\label{bei1}
\sqrt{P_M(\tau^x,\tau^y)}=\tau(|x^*y|)=\tau(\bigl|\,\sqrt{a}\sqrt{c}\,\bigr|)=
\sqrt{P_M(\tau_{a},\tau_{c})}\,.
\end{equation}
Let $s(\tau)$ be the support of $\tau$, that is,
$z=s(\tau)^\perp$ is the central orthoprojection
obeying $I(M,\tau)=Mz$, see \cite[1.10.5.]{Saka:71}. The following is useful in this context\,:
\begin{equation}\label{bei2}
\tau(\bigl|\,\sqrt{a}\sqrt{c}\,\bigr|)=\tau(\sqrt{a}\sqrt{c})\ \Longleftrightarrow\
\sqrt{a}\sqrt{c}\, s(\tau)\geq {\mathbf 0}\,.
\end{equation}
\end{subequations}
\end{example}
In the special case of $M={\mathsf B}({\mathcal H})$ the formula \eqref{bei1} is
known since \cite{Uhlm:76}. With some more effort one even
succeeds in proving \eqref{bei1} without assuming normality. Relating \eqref{bei2},
this is a special case of the fact saying that $\tau(x)=\tau(|x|)$ is equivalent
with  $x\,s(\tau)\geq {\mathbf 0}$, and which can be followed by reasoning about 
polar decomposition of $x$ and invariance of $\tau$ (we omit the details).

\subsubsection{Metric aspects of the Bures distance}\label{met}
From each of the first two items of Theorem \ref{bas5} together with
Definition \ref{budi} and \eqref{bas4}
it gets obvious that $d_{\mathrm{B}}$ is symmetric in its arguments, obeys the triangle inequality,
and is vanishing between
$\nu$ and $\varrho$ if, and only if, $\nu=\varrho$. Thus, $d_{\mathrm{B}}$ yields another
distance function on $M_+^*$. By Theorem \ref{bas5}\,\eqref{bas52}, $d_{\mathrm{B}}$ is the natural
distance function, in each space of $\pi$-fibres.
Also, $d_{\mathrm{B}}$ is topologically
equivalent over bounded subsets to the functional distance $d_1$, which
is defined as $d_1(\nu,\varrho)=\|\nu-\varrho\|_1$, for $\nu,\varrho\in M_+^*$.
For quantitative estimates,
see \cite[{\sc{Proposition 1}},\,formula (1.2)]{AlUh:00.1}.
We do not prove these fact but instead remark that, in case of normal
states/positive linear forms on ${\mathsf W}^*$-algebras, this equivalence results from the
estimates given in
\cite{Bure:69,Arak:71,Arak:72}, essentially. The extension from there
to unital ${\mathsf C}^*$-algebras and their positive linear forms can be easily achieved by means
of the first two items of the following auxiliary result.
\begin{lemma}\label{eqdist}
Let $\nu,\varrho\in M_+^*$, and be $\{\pi,{\mathcal K}\}$ such that the $\pi$-fibres of
$\nu$ and $\varrho$ both exist. Then, the following hold\textup{:}
\begin{enumerate}
\item\label{eqdist1}
$\|\nu-\varrho\|_1=\bigl\|\nu_\pi- \varrho_\pi\bigr\|_1\,;$
\item\label{eqdist2}
$d_{\mathrm{B}}(M|\nu,\varrho)=d_{\mathrm{B}}(\pi(M)^{\,\prime\prime}|\nu_\pi,\varrho_\pi)\,.$
\end{enumerate}
Also, for $\nu,\varrho\in M_+^*$ one has
\begin{enumerate}
\setcounter{enumi}{2}
\item\label{eqdist3}
$P_M(\nu,\varrho)=0\ \Longleftrightarrow\ \nu\perp\varrho\,,$
\end{enumerate}
and for each $\mu,\sigma\in M_+^*$ 
obeying $0\leq\sigma\leq\varrho$, $0\leq\mu\leq\nu$, 
\begin{enumerate}
\setcounter{enumi}{3}
\item\label{eqdist4}
$P_M(\mu,\sigma)+
P_M(\delta\nu,\delta\varrho)\leq P_M(\nu,\varrho)$
\end{enumerate}
is fulfilled, with 
$\delta\nu=\nu-\mu$, $\delta\varrho=\varrho-\sigma$.
\end{lemma}
\begin{proof}
In order to see \eqref{eqdist1} recall first that, since $\pi$ is a $^*$-homomorphism,
according to basic ${\mathsf C}^*$-theory, the unit ball $M_1$ of $M$ is related to the unit ball
$\pi(M)_1$ of $\pi(M)$ through $\pi(M)_1=\pi(M_1)$. From this by the Kaplansky density
theorem
$$\bigl(\pi(M)^{\,\prime\prime}\bigr)_1=
\overline{\pi(M)_1}^{\,\mathrm{st}}=\overline{\pi(M_1)}^{\,\mathrm{st}}\text{ (strong closure) }$$
is seen.
Owing to normality of $\nu_\pi-\varrho_\pi$ we then may conclude as
follows:
\begin{equation*}
\begin{split}
\|\nu-\varrho\|_1 & =\sup_{x\in M,\,\|x\|\leq 1} \bigl|\nu_\pi(\pi(x))-\varrho_\pi(\pi(x))\bigr|\\
& =\sup_{y\in\pi(M)^{\,\prime\prime},\,\|y\|\leq 1} \bigl|\nu_\pi(y)-\varrho_\pi(y)\bigr|= \bigl\|\nu_\pi- \varrho_\pi\bigr\|_1\,.
\end{split}
\end{equation*}
To see \eqref{eqdist2}, apply Theorem \ref{bas5}\,\eqref{bas51} to $\nu_\pi$ and $\varrho_\pi$ in respect of the identity representation
of the $vN$-algebra
$\pi(M)^{\,\prime\prime}$. Under these premises
Theorem \ref{bas5}\,\eqref{bas51} yields
$d_{\mathrm{B}}(\pi(M)^{\,\prime\prime}|\nu_\pi,\varrho_\pi)=\sup_{u\in {\mathcal U}(\pi(M)^{\,\prime})}\|u\psi-\varphi\|$, for each $\psi\in{\mathcal S}_{\mathrm{id},\pi(M)^{\,\prime\prime}}(\varrho_\pi)$ and $\varphi\in
{\mathcal S}_{\mathrm{id},\pi(M)^{\,\prime\prime}}(\nu_\pi)$.
In view of Remark \ref{rem1}\,\eqref{rem11} this is $d_{\mathrm{B}}(M|\nu,\varrho)$ as asserted by
Theorem \ref{bas5}\,\eqref{bas51} when applied with respect to $M$, $\pi$, $\nu$
and $\varrho$.

Relating \eqref{eqdist4}, note that for each invertible $x\in M_+$ obviously
$\nu(x)\varrho(x^{-1})=\{\nu-\mu\}(x)\{\varrho-\sigma\}(x^{-1})+\nu(x)\sigma(x^{-1})+
\mu(x)\varrho(x^{-1})+\mu(x)\sigma(x^{-1})$ holds.
From this by Lemma \ref{bas3} and positivity of all terms the assertion follows.

To see \eqref{eqdist3}, note that $P_M(\nu,\varrho)=0$ is equivalent to  $h_{\psi,\varphi}^\pi=0$, by
Lemma \ref{bas3}, which is the same as $\langle x\psi,y\varphi\rangle=0$, for all
$x,y\in \pi(M)^{\,\prime}$. Hence, $P_M(\nu,\varrho)=0$ is equivalent to  $p_\pi(\psi)\perp
p_\pi(\varphi)$. Therefore, since $p_\pi(\psi)=s(\varrho_\pi)$ and
$p_\pi(\varphi)=s(\nu_\pi)$ holds, with the
support orthoprojections $s(\nu_\pi)$ and $s(\varrho_\pi)$ of the normal positive linear forms
$\nu_\pi$ and $\varrho_\pi$ (see Remark \ref{rem1}\,\eqref{rem11}),
$P_M(\nu,\varrho)=0$ is equivalent to orthogonality of $\nu_\pi$ with $\varrho_\pi$.
The latter is the same as
$\|\nu_\pi-\varrho_\pi\|_1=\|\nu_\pi\|_1+\|\varrho_\pi\|_1$.
In view of \eqref{eqdist1}
equivalence of $P_M(\nu,\varrho)=0$ to 
$\|\nu-\varrho\|_1=\|\nu\|_1+\|\varrho\|_1$ follows.
\end{proof}
Relating
Theorem \ref{bas5}\,\eqref{bas53}, a remarkable fact firstly acknowledged in \cite{Arak:72}
says that $d_{\mathrm{B}}(M|\nu,\varrho)=\|\psi-\varphi\|$
can happen, for each $\nu,\varrho\in M_+^*$, see also \cite[Appendix 7.]{Albe:92.1}.
\begin{theorem}\label{positiv1}
Suppose $\{\pi,{\mathcal K}\}$ is such that the $\pi$-fibres of
$\nu,\varrho\in M_+^*$ both exist. There are $\varphi_0\in {\mathcal S}_{\pi,M}(\nu)$ and
$\psi_0\in {\mathcal S}_{\pi,M}(\varrho)$ with $d_{\mathrm{B}}(M|\nu,\varrho)=\bigl\|\psi_0-\varphi_0\bigr\|$.
\end{theorem}
\begin{proof}
For $\varphi\in {\mathcal S}_M(\nu)$, $\psi\in {\mathcal S}_M(\varrho)$ consider
$h_{\psi,\varphi}^\pi$ as defined in \eqref{hform}. Then
$\bigl\|h_{\psi,\varphi}^\pi\bigr\|_1=\bigl\|\bigl|h_{\psi,\varphi}^\pi\bigr|\bigr\|_1$ holds,
see \eqref{sub}. Hence, in the notations from above
$$\bigl\|h_{\psi,\varphi}^\pi\bigr\|_1=\bigl\langle v_{\psi,\varphi}^{\pi*}\psi,\varphi\bigr
\rangle=h_{\psi,\varphi}^\pi\bigl(v_{\psi,\varphi}^{\pi*}\bigr).$$

It suffices to show that partial
isometries $v,w\in \pi(M)^{\,\prime}$ exist which obey
$v^*v\geq p_\pi^\prime(\varphi)$, $w^*w\geq p_\pi^\prime(\psi)$ and $h_{\psi,\varphi}^\pi\bigl(v_{\psi,\varphi}^{\pi*}\bigr)=
h_{\psi,\varphi}^\pi(v^*w)$. In fact, in view of \eqref{bas4}, for
$\psi_0=w\psi$ and $\varphi_0=v\varphi$
one then has $\psi\in {\mathcal S}_{\pi,M}(\varrho)$ and $\varphi\in {\mathcal S}_{\pi,M}(\nu)$.
By the above and in twice applying Lemma \ref{bas3}, we may conclude as follows\,:
$$h_{\psi_0,\varphi_0}^\pi({\mathbf 1})=h_{\psi,\varphi}^\pi(v^*w)=h_{\psi,\varphi}^\pi\bigl(v_{\psi,\varphi}^{\pi*}\bigr)=
\bigl\|h_{\psi,\varphi}^\pi\bigr\|_1=\sqrt{P_M(\nu,\varrho)}=\bigl\|h_{\psi_0,\varphi_0}^\pi\bigr\|_1\,.$$
Hence $h_{\psi_0,\varphi_0}^\pi({\mathbf 1})=
\bigl\|h_{\psi_0,\varphi_0}^\pi\bigr\|_1$, from which
$h_{\psi_0,\varphi_0}^\pi\geq 0$ follows. In view of
Theorem \ref{bas5}\,\eqref{bas53}, the constructed $\varphi_0$ and $\psi_0$
can be taken to meet our demands.

Recall that for a $vN$-algebra $N$ on a Hilbert space ${\mathcal H}$
there is a largest central orthoprojection
$z\in N\cap N^{\,\prime}$ such that
the $vN$-algebra $(N^{\,\prime})z$ over $z{\mathcal H}$ is finite.
Thus, in case of $z\not={\mathbf 1}$,
$(N^{\,\prime}){z^\perp}$ is properly infinite over $z^\perp{\mathcal H}$.
Moreover, on $c{\mathcal K}$ one has
$$(N^{\,\prime})c=\bigl((N^{\,\prime})c\bigr)^{\,\prime}\,,$$
for each central orthoprojection $c$
(the outer commutant
refers to $c{\mathcal K}$).
Applying this to $N=\pi(M)^{\,\prime}$ with $c=z$ or $c=z^\perp$, and having in mind
Remark \ref{rem1}\,\eqref{rem11}, we see that we can content ourselves  with constructing the partial
isometries in question only for
$vN$-algebras with either finite or properly infinite commutants.

In line with this, let $M$ be the $vN$-algebra in question, with $M^{\,\prime}$ either finite or
properly infinite, and let $\nu$ and $\varrho$ be implemented by vectors $\varphi$ and $\psi$.
We then will make use of the notational conventions of
Remark \ref{rem1}\,\eqref{rem12}, and are going to construct $v$ and $w$ under these premises now.

Suppose $M^{\,\prime}$ to be finite.
By \eqref{sub3} we have
$s\bigl(\bigl|h_{\psi,\varphi}^*\bigr|\bigr)\sim
s\bigl(\bigl|h_{\psi,\varphi}\bigr|\bigr)$.
In a finite $vN$-algebra, by another standard fact, see \cite[2.4.2.]{Saka:71}, this condition implies
the following to hold, in any case:
\begin{equation}\label{dec}
s\bigl(\bigl|h_{\psi,\varphi}^*\bigr|\bigr)^\perp
\sim s\bigl(\bigl|h_{\psi,\varphi}\bigr|\bigr)^\perp\,.
\end{equation}
Hence, there is $m\in M^{\,\prime}$ obeying $m^*m=s\bigl(\bigl|h_{\psi,\varphi}\bigr|\bigr)^\perp$
and  $m m^*=s\bigl(\bigl|h_{\psi,\varphi}^*\bigr|\bigr)^\perp$. Define $w=v_{\psi,\varphi}^*+m^*$.
Then, $w$ is unitary, $w\in {\mathcal U}(M^{\,\prime})$, with
$w s\bigl(\bigl|h_{\psi,\varphi}^*\bigr|\bigr)=v_{\psi,\varphi}^*$. Since by
polar decomposition $h_{\psi,\varphi}(w)=h_{\psi,\varphi}\bigl(v_{\psi,\varphi}^*\bigr)$
is fulfilled, $v={\mathbf 1}$ can be chosen.

Suppose a properly infinite $M^{\,\prime}$. Then
$p\in M^{\,\prime}$ with $p\sim p^\perp\sim {\mathbf 1}$ exists.
By \eqref{sub1},
$p^\prime(\varphi)-s\bigl(\bigl|h_{\psi,\varphi}\bigr|\bigr)\in M^{\,\prime}$ is an
orthoprojection. Hence
$\bigl\{p^\prime(\varphi)-s\bigl(\bigl|h_{\psi,\varphi}\bigr|\bigr)\bigr\}
\prec p$ and $p'(\psi)\prec p^\perp$ (`$\prec$' be the
Murray-von Neumann comparability relation). Hence,
there are
$v_1,v_2\in M^{\,\prime}$ with
$v_1^*v_1=\bigl\{p^\prime(\varphi)-s\bigl(\bigl|h_{\psi,\varphi}\bigr|\bigr)\bigr\}$,
$v_1v_1^*\leq p$, $v_2^*v_2=p'(\psi)$, $v_2 v_2^*\leq p^\perp$. Define
$w=v_2$, $v=v_1+v_2 v_{\psi,\varphi}$. Since $v_1^*v_2=v_2^*v_1={\mathbf 0}$ holds, the
above and \eqref{sub2}, \eqref{sub1} imply
$v^*v=v_1^*v_1+v_{\psi,\varphi}^*v_2^*v_2v_{\psi,\varphi}=
\bigl\{p'(\varphi)-s\bigl(\bigl|h_{\psi,\varphi}\bigr|\bigr)\bigr\}+
s\bigl(\bigl|h_{\psi,\varphi}\bigr|\bigr)=p'(\varphi)$ and $w^*w=p'(\psi)$.
With the help of \eqref{sub2} again, we finally see that $v^*w=
v_{\psi,\varphi}^*v_2^*v_2=v_{\psi,\varphi}^*p'(\psi)=
v_{\psi,\varphi}^*s\bigl(\bigl|h_{\psi,\varphi}^*\bigr|\bigr)=
v_{\psi,\varphi}^*$ holds.
\end{proof}
\begin{remark}\label{subadd}
From Lemma \ref{eqdist}\,\eqref{eqdist4} the convexity/subadditivity properties
of $P$ are obvious: $P$ is subadditive. On the other hand, by the 
scaling behavior of $P$ on the cone $M_+^*$ together with Theorem \ref{bas5}\,\eqref{bas52} and 
formula \eqref{pcont.1}, it is easily inferred that $d_{\mathrm{B}}^{\,2}$ (resp.~$\sqrt{P}$)  
is jointly convex (resp.~jointly concave), on each convex subset of $M_+^*$,  
see also in \cite{AlUh:84}, and \cite[Remark 2\,(3)]{AlUh:00.1} for
some additional superadditivity property of $\sqrt{P}$.
\end{remark}
\subsubsection{Minimal pairs of positive linear forms}\label{minpairs}
Let $\varphi_0$ and $\psi_0$ in the $\pi$-fibres of $\nu,\varrho\in M_+^*$
be chosen as to satisfy the hypothesis of Theorem \ref{positiv1}.
Define $\varrho^\perp\in M_+^*$ as follows:
\begin{equation}\label{deforth}
\eta\in {\mathcal S}_{\pi,M}(\varrho^\perp),\text{ with }\eta=s(h_{\psi_0,\varphi_0}^\pi)^\perp\psi_0\,.
\end{equation}
Owing to $s(h_{\psi_0,\varphi_0}^\pi)^\perp\in \pi(M)^{\,\prime}$ and
$\psi_0\in {\mathcal S}_{\pi,M}(\varrho)$, $0\leq \varrho^\perp\leq \varrho$.
Also $h_{\eta,\varphi_0}^\pi=0$, since
$h_{\eta,\varphi_0}^\pi(z)=h_{\psi_0,\varphi_0}^\pi(zs(h_{\psi_0,\varphi_0}^\pi)^\perp)
=0$,
for each $z\in \pi(M)^{\,\prime}$. By Lemma \ref{bas3}
$P_M(\nu,\varrho^\perp)=0$, and thus
$\varrho^\perp$ and $\nu$ are mutually orthogonal, by Lemma \ref{eqdist}\,\eqref{eqdist3}.

On the other hand, let $\sigma\in M_+^*$ with $\sigma\leq \varrho$ be orthogonal to $\nu$.
Then, by standard facts there exists $z\in \pi(M)^{\,\prime}$ with $\|z\|\leq 1$ and
$z\psi_0\in {\mathcal S}_{\pi,M}(\sigma)$. Also, by
Lemma \ref{eqdist}\,\eqref{eqdist3}, $P_M(\nu,\sigma)=0$. Hence, for $\eta'=z\psi_0$,
$h_{\eta',\varphi_0}^\pi=0$ by Lemma \ref{bas3}.
Thus in particular $h_{\psi_0,\varphi_0}^\pi(z^*z)=0$.
By positivity of $h_{\psi_0,\varphi_0}^\pi$ and $\|z\|\leq 1$
one concludes $z^*z\leq {\mathbf 1}-s(h_{\psi_0,\varphi_0}^\pi)$.
In view of \eqref{deforth} then
$\sigma\leq \varrho^\perp$ follows. Thus we arrive at the following result,
which goes back to \cite{Arak:72},
see \cite[\S 2]{Albe:92.1} and \cite{Hein:91} for applications.
\begin{theorem}\label{ortho}
For each given pair $\{\nu,\varrho\}$ of positive linear
forms on a unital ${\mathsf C}^*$-algebra, among all positive linear forms $\sigma$
with $\sigma\perp \nu$
and $\sigma\leq \varrho$ there exists
a largest element, $\sigma=\varrho^\perp$.
\end{theorem}
For given pair $\{\nu,\varrho\}$ there has to exist also
a largest positive linear form subordinate to $\nu$ and orthogonal to $\varrho$, which we call $\nu^\perp$.
Thus, for each pair $\{\nu,\varrho\}$ of positive linear forms
both $\nu^\perp$ and $\varrho^\perp$ have an invariant meaning. The
derivation of the previous result then also shows that, provided the $\pi$-fibres of $\nu$ and $\varrho$ both exist, with respect to the pair
$\{\nu_\pi,\varrho_\pi\}$ over $\pi(M)^{\,\prime\prime}$ the following must be valid:
\begin{equation}\label{perp}
(\varrho_\pi)^\perp=(\varrho^\perp)_\pi,\,(\nu_\pi)^\perp=(\nu^\perp)_\pi\,.
\end{equation}
On the other hand, since
$\varrho-\varrho^\perp\geq 0$ holds, we may also consider the pair $\{\nu, \varrho-\varrho^\perp\}$, and
we may ask for $(\varrho-\varrho^\perp)^{\perp'}$ and $\nu^{\perp'}$, where ${\perp'}$ now
refers to the $\perp$-operation with respect to the pair $\{\nu,\varrho-\varrho^\perp\}$.
\begin{lemma}\label{perp1}
For each pair $\{\nu,\varrho\}$ the following facts hold\textup{:}
\begin{enumerate}
\item\label{perp1.1}
$(\varrho-\varrho^\perp)^{\perp'}=0,\ \nu^{\perp'}=\nu^\perp\,;$
\item\label{perp1.3}
$d_{\mathrm{B}}(\nu,\varrho)^2= d_{\mathrm{B}}(\nu-\nu^\perp,\varrho-\varrho^\perp)^2 +
d_{\mathrm{B}}(\nu^\perp,\varrho^\perp)^2\,;$
\item\label{perp1.4}
$P_M(\nu,\varrho)=P_M(\nu-\nu^\perp,\varrho-\varrho^\perp)\,.$
\end{enumerate}
Especially, the pair $\{\nu-\nu^\perp,\varrho-\varrho^\perp\}$ is the least element of all pairs
$\{\theta,\theta'\}$ of positive linear forms with $\theta\leq \nu$ and
$\theta'\leq \varrho$, and which obey $P_M(\theta,\theta')=P_M(\nu,\varrho)$.
\end{lemma}
\begin{proof}
By \eqref{perp}, Lemma \ref{eqdist} and Theorem \ref{ortho} we may content ourselves  with proving the assertions
for a $vN$-algebra $M$ and normal positive linear forms $\nu$ and $\varrho$ with nontrivial
fibres. For normal forms orthogonality simply means orthogonality of the respective supports.
Hence, $0\leq \sigma\leq \varrho-\varrho^\perp$ and $\sigma\perp \nu$
imply $\sigma+\varrho^\perp\leq \varrho$ and $\sigma+\varrho^\perp\perp \nu$. By maximality of $\varrho^\perp$
then $\sigma=0$, and thus the first of the relations of \eqref{perp1.1} is seen.
Note that $\nu^{\perp'}\leq \nu$, with $\nu^{\perp'}\perp \varrho-\varrho^\perp$.
From the former owing to $\varrho^\perp\perp\nu$ then $\varrho^\perp\perp\nu^{\perp'}$ is seen,
and thus the two orthogonality relations may be summarized into
$\nu^{\perp'}\perp (\varrho-\varrho^\perp)+\varrho^\perp=\varrho$. Hence, in view of
Theorem \ref{ortho} then $\nu^{\perp'}\leq \nu^\perp$ follows. On the other hand, by definition,
$\nu^{\perp}$ is subordinate to $\nu$ and must also be orthogonal to $\varrho-\varrho^\perp\leq \varrho$. Hence,
$\nu^{\perp'}\geq \nu^\perp$ by Theorem \ref{ortho} when applied to
$\{\nu, \varrho-\varrho^\perp\}$. In summarizing the second
relation of \eqref{perp1.1} follows.

Let $\psi_0\in {\mathcal S}_M(\varrho)$ and $\varphi_0\in {\mathcal S}_M(\nu)$ be chosen
as in Theorem \ref{positiv1}. By Theorem \ref{bas5}\,\eqref{bas53} and
\eqref{deforth} one has
$s^\perp\psi_0\in {\mathcal S}_M(\varrho^\perp)$, with $s=s(h_{\psi_0,\varphi_0})\in M^{\,\prime}$.
Note that $h_{\psi_0,\varphi_0}=h_{\psi_0,\varphi_0}((\cdot)s)=h_{s\psi_0,\varphi_0}$. By
Theorem \ref{bas5}\,\eqref{bas53} $d_{\mathrm{B}}(\nu,\varrho-\varrho^\perp)=\|\varphi_0-s\psi_0\|$
follows. By orthogonality of $\varrho^\perp$ with $\nu$ we especially must have
$s^\perp\psi_0\perp \varphi_0$. Hence also $s^\perp\psi_0\perp (\varphi_0-s\psi_0)$, and thus we obtain
$$d_{\mathrm{B}}(\nu,\varrho)^2=\|\varphi_0-s\psi_0-s^\perp\psi_0\|^2=d_{\mathrm{B}}(\nu,\varrho-\varrho^\perp)^2+
\|s^\perp\psi_0\|^2=d_{\mathrm{B}}(\nu,\varrho-\varrho^\perp)^2+\|\varrho^\perp\|_1\,.$$
We may apply this
to $\{\varrho-\varrho^\perp,\nu\}$ accordingly, with the result
$d_{\mathrm{B}}(\varrho-\varrho^\perp,\nu)^2=d_{\mathrm{B}}(\varrho-\varrho^\perp,\nu-\nu^{\perp'})^2+
\|\nu^{\perp'}\|_1$.
By \eqref{perp1.1} and by symmetry of $d_{\mathrm{B}}$,
substitution of the latter into the former relation yields
$d_{\mathrm{B}}(\nu,\varrho)^2= d_{\mathrm{B}}(\nu-\nu^\perp,\varrho-\varrho^\perp)^2 +\|\nu^\perp\|_1+\|\varrho^\perp\|_1$.
Since $\varrho^\perp\perp\nu^\perp$ and
Lemma \ref{eqdist}\,\eqref{eqdist3} and formula \eqref{budi} imply
$\|\nu^\perp-\varrho^\perp\|_1=
\|\nu^\perp\|_1+\|\varrho^\perp\|_1=d_{\mathrm{B}}(\nu^\perp,\varrho^\perp)^2$,
from this \eqref{perp1.3} follows. In view of \eqref{budi},
\eqref{perp1.4} is seen.
Finally, by Lemma \ref{eqdist}\,\eqref{eqdist4},
$P_M(\theta,\theta')\leq P_M(\nu,\theta')\leq P_M(\nu,\varrho)$ as well as
$P_M(\theta,\theta')\leq P_M(\nu,\theta')+P_M(\nu,\varrho-\theta')\leq
P_M(\nu,\varrho)$ must be fulfilled. Therefore,
$P_M(\theta,\theta')=P_M(\nu,\varrho)$ implies $P_M(\nu,\varrho-\theta')=0$.
By  Lemma \ref{eqdist}\,\eqref{eqdist3} then $\nu \perp (\varrho-\theta')$ follows.
Since $0\leq (\varrho-\theta')\leq \varrho$ holds, in view of Theorem \ref{ortho}
$\varrho^\perp\geq (\varrho-\theta')$ follows. Hence, $\theta'\geq \varrho-\varrho^\perp$.
Proceeding analogously for $\theta$ instead of $\theta'$ yields $\theta\geq \nu-\nu^\perp$.
Thus also the final assertion is true.
\end{proof}
\begin{definition}\label{minpair}
Refer to $\{\omega,\sigma\}$ as {\em minimal pair} if
$\omega^\perp=\sigma^\perp=0$ holds, and refer to $\{\nu-\nu^\perp,\varrho-\varrho^\perp\}$ as
$\{\nu,\varrho\}$-associated minimal pair.
\end{definition}
\begin{remark}\label{remminpair}
It is obvious from the last part of Lemma \ref{perp1} that the
$\{\nu,\varrho\}$-associated minimal pair is the unique minimal pair
$\{\theta,\theta'\}$ obeying $\theta\leq \nu$,
$\theta'\leq \varrho$ and $P_M(\theta,\theta')=P_M(\nu,\varrho)$, see \cite[Lemma 2.1]{Albe:92.1}.
\end{remark}
There is yet another remarkable uniqueness result in context of these structures\,:
\begin{theorem}{\cite[Theorem 3.1]{Albe:92.1}}\label{rem23}
For given pair $\{\nu,\varrho\}$ of positive linear forms on a unital
${\mathsf C}^*$-algebra $M$, there
is a unique $g\in M^*$ obeying
\begin{enumerate}
\item\label{a}
$g({\mathbf 1})=\sqrt{P_M(\nu,\varrho)}\,,$
\item\label{b}
$|\,g(y^*x)|^2\leq \{\nu-\nu^\perp\}(y^*y)\{\varrho-\varrho^\perp\}(x^*x)\,,\ x,y\in M\,.$
\end{enumerate}
\end{theorem}
\section{Forthcoming results}\label{forth}
\subsection{Setting the problems}\label{forth1}
The problems of this paper naturally came along when analyzing in context of
the hypothesis of Theorem \ref{positiv1} an example recently noticed in \cite{Pelt:00.1}
and some other optimization problem, see
\cite[Theorem 6\,(2)]{AlUh:00.1}.
\begin{problem}\label{prob}
\begin{enumerate}
\item\label{prob1}
Despite the given proof of Theorem \ref{positiv1}, we are missing
a true criterion along which one can decide whether or not, for given individual
vector $\varphi_0$ in the $\pi$-fibre of $\nu$, the condition in the
hypothesis of Theorem \ref{positiv1} could be satisfied at all, by some vector $\psi_0$ in the
$\pi$-fibre of $\varrho$.
\item\label{prob2}
Let ${\mathcal S}_{\pi,M}(\nu|\varrho)$ be the subset of all $\varphi_0$ in the $\pi$-fibre of $\nu$
for which the question in \eqref{prob1} can be affirmatively answered. Then, provided $\varphi_0\in
{\mathcal S}_{\pi,M}(\nu|\varrho)$ is fulfilled, is there a formula giving a particular
$\psi_0$ to that $\varphi_0$\,?
\item\label{prob3}
Especially, by Theorem \ref{positiv1} the question is left open whether
${\mathcal S}_{\pi,M}(\nu|\varrho)={\mathcal S}_{\pi,M}(\nu)$ might happen, on a given unital
${\mathsf C}^*$-algebra.
\item\label{prob4}
And finally, provided for some $M$ the answer in \eqref{prob3} can be in the negative,
what about the structure of those pairs $\{\nu,\varrho\}$ with
${\mathcal S}_{\pi,M}(\nu|\varrho)\not={\mathcal S}_{\pi,M}(\nu)$\,?
\end{enumerate}
\end{problem}
Thereby, in line with Remark \ref{rem1}\,\eqref{rem11} and
Lemma \ref{eqdist}\,\eqref{eqdist2}, it is reasonable to content ourselves with considering these
questions for $vN$-algebras with normal positive linear forms which can be implemented by vectors.
In making reference to some
ideas about hereditarity (we omit all these details), we even can be assured that
without essential loss of generality we might confine ourselves
to considering $d_{\mathrm{B}}$
over normal positive linear forms over a standard form $vN$-algebra $M$.
Therefore, the problems in questions will be dealt with
and answered examplarily at least for such cases. Then also the
notational conventions of Remark \ref{rem1}\,\eqref{rem12} tacitly will be in use, and
${\mathcal S}_M(\nu|\varrho)$ instead of ${\mathcal S}_{\pi,M}(\nu|\varrho)$ will be written.
\begin{remark}\label{fin}
On a {\em finite} standard form $vN$-algebra $M$,
problems \eqref{prob1},\,\eqref{prob3} and
\eqref{prob4} trivialize, for then always
\begin{equation}\label{vorbild}
{\mathcal S}_M(\nu)={\mathcal S}_M(\nu|\varrho)=\bigl\{u\varphi_0:\,u^*u={\mathbf 1},\,u\in M^{\,\prime}\bigr\}
\end{equation}
holds, with arbitrary $\varphi_0\in {\mathcal S}_M(\nu)$. In fact, for finite $M$, $M^{\,\prime}$ is finite, too.
Since in a finite $vN$-algebra each partial isometry admits unitary extensions,
for each orthoprojection $q\in M^{\,\prime}$ one has
$\{v\in M^{\,\prime}:v^*v=q \}={\mathcal U}(M^{\,\prime})q$. Especially then also ${\mathcal U}(M^{\,\prime})=\{v\in M^{\,\prime}:v^*v={\mathbf 1} \}$ follows. Hence, in view of \eqref{bas4}
each fibre is the ${\mathcal U}(M^{\,\prime})$-orbit of each of its vectors. From this and
Theorem \ref{positiv1} then \eqref{vorbild} follows.
\end{remark}
\subsection{The general criterion}\label{forth2}
The first two problems of Problem \ref{prob} can be answered
under the rather general settings of a unital ${\mathsf C}^*$-algebra $M$, a pair
$\{\nu,\varrho\}$ of positive linear forms, and $\varrho^\perp$ as determined by
Theorem \ref{ortho}. Under the same premises as in Theorem \ref{positiv1}, we
have the following result.
\begin{theorem}\label{sesq*}
A vector $\varphi_0\in {\mathcal
S}_{\pi,M}(\nu)$ belongs to ${\mathcal
S}_{\pi,M}(\nu|\varrho)$ if, and only if,
there are $\eta\in{\mathcal S}_{\pi,M}(\varrho^\perp)$ and
$\psi\in {\mathcal S}_{\pi,M}(\varrho)$ satisfying
\begin{subequations}\label{crit}
\begin{equation}\label{positiv1*}
p_\pi^{\,\prime}(\eta)\leq s(|h_{\psi,\varphi_0}^\pi|)^\perp\,.
\end{equation}
If this condition is fulfilled, then the special vector
\begin{equation}\label{impmod.1*}
\psi_0=v_{\psi,\varphi_0}^{\pi*}\psi+\eta
\end{equation}
belongs to ${\mathcal S}_{\pi,M}(\varrho)$, that is,
$d_{\mathrm{B}}(\nu,\varrho)=\|\psi_0-\varphi_0\|$ is satisfied.
Thereby, for given $\varphi_0$ and $\eta$, \eqref{positiv1*} remains
true for any $\psi\in {\mathcal
S}_{\pi,M}(\varrho)$, and
$\psi_0$ as given by formula \eqref{impmod.1*} is the same, for each such $\psi$.
\end{subequations}
\end{theorem}
For the proof, it is useful to take notice of the following auxiliary results first.
\begin{lemma}\label{sesq}
Let $\nu,\varrho\in M_{*+}$ and $\varphi\in {\mathcal S}_{\pi,M}(\nu)$,
$\psi\in {\mathcal S}_{\pi,M}(\varrho)$.
The following hold\textup{:}
\begin{enumerate}
\item\label{positiv-2}
$p_\pi'\bigl(v_{\psi,\varphi}^{\pi*}\psi\bigr)=s\bigl(\bigl|h_{\psi,\varphi}^\pi
\bigr|\bigr)\,;$
\item\label{positiv-1}
$\forall\,\tilde{\psi}\in {\mathcal S}_{\pi,M}(\varrho)\,:\
\bigl|h_{\tilde{\psi},\varphi}^\pi\bigr|=\bigl|h_{\psi,\varphi}^\pi\bigr|,\,
v_{\tilde{\psi},\varphi}^\pi=wv_{\psi,\varphi}^\pi,\ \text{with}\ w\in
\pi(M)^{\,\prime}\ \text{obeying}\\
w^*w=p_\pi'(\psi),\,\tilde{\psi}=w\psi\,;$
\item\label{positiv0*}
$v_{\psi,\varphi}^{\pi*}\psi,\,s\bigl(|h_{\psi,\varphi}^{\pi*}|\bigr)\psi\in
{\mathcal S}_{\pi,M}(\varrho-\varrho^\perp),\
s\bigl(|h_{\psi,\varphi}^{\pi*}|\bigr)^\perp\psi\in
{\mathcal S}_{\pi,M}(\varrho^\perp)\,.$
\end{enumerate}
\end{lemma}
\begin{proof}
Note that $s\bigl(\bigl|h_{\psi,\varphi}^\pi\bigr|\bigr)\leq
p_\pi'\bigl(v_{\psi,\varphi}^{\pi*}\psi\bigr)$, by \eqref{sub0.5}, and
$p_\pi'\bigl(v_{\psi,\varphi}^{\pi*}\psi\bigr)\leq
s\bigl(\bigl|h_{\psi,\varphi}^\pi\bigr|\bigr)$ by \eqref{sub1}. Hence
\eqref{positiv-2} follows. Relating \eqref{positiv-1}, remark that by \eqref{bas4} $w\in
\pi(M)^{\,\prime}$ obeying $w^*w=p_\pi'(\psi)$ and $\tilde{\psi}=w\psi$ exists
and is unique. In view of this and \eqref{sub2},
$u=wv_{\psi,\varphi}^\pi$ is a partial isometry of $\pi(M)^{\,\prime}$, with
$u^*u=v_{\psi,\varphi}^{\pi*}v_{\psi,\varphi}^\pi=
s\bigl(\bigl|h_{\psi,\varphi}^\pi\bigr|\bigr)$.  From this and
$h_{\tilde{\psi},\varphi}^\pi=h_{\psi,\varphi}^\pi
((\cdot)w)=\bigl|h_{\psi,\varphi}^\pi\bigr|
( (\cdot)u)$ by uniqueness of
the polar decomposition the validity of \eqref{positiv-1} follows.

Let $\psi_0,\varphi_0$ be as in Theorem \ref{positiv1}.
By Theorem \ref{bas}\,\eqref{bas53}, $h_{\psi_0,\varphi_0}^\pi\geq 0$, and according to
\eqref{deforth}, $s\bigl(h_{\psi_0,\varphi_0}^\pi\bigr)\psi_0\in{\mathcal
S}_{\pi,M}(\varrho-\varrho^\perp)$. By
\eqref{positiv-1} one has $w,u\in \pi(M)^{\,\prime}$,
with $w^*w=p_\pi'(\psi_0)$, $\psi=w\psi_0$, $u^*u=p_\pi'(\varphi_0)$,
$\varphi=u\varphi_0$ and obeying $v_{\psi,\varphi_0}^\pi=w
v_{\psi_0,\varphi_0}^\pi$, $v_{\varphi,\psi}^\pi=u v_{\varphi_0,\psi}^\pi$.
The partial isometry $v_{\alpha,\beta}$ of the polar
decomposition of $h_{\alpha,\beta}$ has to obey
$v_{\alpha,\beta}^*=v_{\beta,\alpha}$. Hence,
$v_{\psi,\varphi}^{\pi*}=v_{\varphi,\psi}^\pi=u
v_{\varphi_0,\psi}^\pi=u v_{\psi,\varphi_0}^\pi=
u v_{\psi_0,\varphi_0}^{\pi*} w^*$ follows. From this together with the
other properties of $w,u$ one gets
\begin{equation}\label{star}
v_{\psi,\varphi}^{\pi*}\psi=
u v_{\psi_0,\varphi_0}^{\pi*} \psi_0\,.
\end{equation}
By $h_{\psi_0,\varphi_0}^\pi\geq 0$ one has
$v_{\psi_0,\varphi_0}^\pi=s\bigl(h_{\psi_0,\varphi_0}^\pi\bigr)$. Thus, from
$u^*u=p_\pi'(\varphi_0)\geq
s\bigl(h_{\psi_0,\varphi_0}^\pi\bigr)$ in view of
$s\bigl(h_{\psi_0,\varphi_0}^\pi\bigr)\psi_0\in{\mathcal
S}_M(\varrho-\varrho^\perp)$ and \eqref{star} also
$v_{\psi,\varphi}^{\pi*}\psi  \in{\mathcal
S}_{\pi,M}(\varrho-\varrho^\perp)$ follows. In view of \eqref{sub2} and \eqref{bas4}
then $s\bigl(|h_{\psi,\varphi}^{\pi*}|\bigr)\psi\in{\mathcal
S}_{\pi,M}(\varrho-\varrho^\perp)$ is obtained. From this \eqref{positiv0*} follows.
\end{proof}
\begin{proof}{(of Theorem \ref{sesq*})}
Suppose $\eta\in {\mathcal S}_{\pi,M}(\varrho^\perp)$ satisfies \eqref{positiv1*}
for some $\psi\in {\mathcal S}_M(\varrho)$,
and given $\varphi_0\in {\mathcal S}_{\pi,M}(\nu)$. By
Lemma \ref{sesq}\,\eqref{positiv-1} the same premises then hold
with respect to any other vector $\psi$ in the $\pi$-fibre of $\varrho$. Also,
by Lemma \ref{sesq}\,\eqref{positiv-2}, one has
$p_\pi^\prime(v_{\psi,\varphi_0}^{\pi*} \psi)=s\bigl(|h_{\psi,\varphi_0}^\pi|\bigr)$. From
this and \eqref{positiv1*} in view of Lemma \ref{sesq}\,\eqref{positiv0*} then $\psi_0\in {\mathcal
S}_{\pi,M}(\varrho)$ follows. By
Theorem \ref{ortho}, $\varrho^\perp \perp \nu$. This is the same as
$p_\pi(\eta)\perp p_\pi(\varphi_0)$, and thus in view of \eqref{sub0.5} and by construction
of $\psi_0$, in accordance with \eqref{impmod.1*}
$h_{\psi_0,\varphi_0}^\pi=h_{\psi,\varphi_0}^\pi\bigl((\cdot)
v_{\psi,\varphi_0}^{\pi*}\bigr)=|h_{\psi,\varphi_0}^\pi|\geq 0$.
By Theorem \ref{bas5}\,\eqref{bas53} we see $d_{\mathrm{B}}(\nu,\varrho)=\|\psi_0-\varphi_0\|$, for
$\psi_0$ of \eqref{impmod.1*}.
Note that our situation with $\varphi_0,\,\psi_0$ and $\psi$ can be easily adapted to a context where
relation \eqref{star} can be used. In fact, choosing $\varphi=\varphi_0$ there,
in the notation of the previous proof
one has $u=p_\pi^\prime(\varphi_0)$, and in view of
$p_\pi^\prime(v_{\psi,\varphi_0}^{\pi*} \psi)=s\bigl(|h_{\psi,\varphi_0}^\pi|\bigr)\leq p_\pi^\prime(\varphi_0)$
(see above and \eqref{sub1}) the relation  \eqref{star}
yields $v_{\psi,\varphi_0}^{\pi*} \psi=v_{\psi_0,\varphi_0}^{\pi*} \psi_0$, for each other
$\psi\in {\mathcal
S}_{\pi,M}(\varrho)$. Hence, $\psi_0$ of
\eqref{impmod.1*} does not depend on the special $\psi$, provided
$\varphi_0$, $\eta$ are held fixed.

To see that \eqref{positiv1*} is necessary, let
$\psi_0,\varphi_0$ obey $d_{\mathrm{B}}(\nu,\varrho)=\|\psi_0-\varphi_0\|$.
By Theorem \ref{bas5}\,\eqref{bas53},
$h_{\psi_0,\varphi_0}^\pi\geq 0$, and thus
$s(h_{\psi_0,\varphi_0}^\pi)^\perp\psi_0\in {\mathcal
S}_{\pi,M}(\varrho^\perp)$, by Lemma \ref{sesq}\,\eqref{positiv0*}.
Obviously, for $\psi=\psi_0$ the condition \eqref{positiv1*}
can be satisfied with $\eta=s(h_{\psi_0,\varphi_0}^\pi)^\perp\psi_0$. In view of the
above this completes the proof.
\end{proof}

\subsection{Specific applications}\label{forth3}
In the following
$M$ will be a $vN$-algebra acting in
standard form on a
Hilbert space ${\mathcal
H}$, with cyclic and
separating vector ${\varOmega}\in {\mathcal H}$.
The previous results then make sense with respect to the identity representation, and can be considered
for each pair $\{\nu,\varrho\}$, $\nu,\varrho\in M_{*+}$ (normal positive linear forms).

\subsubsection{Modular cone and implementing vectors}\label{covN}
Recall a few basic facts from modular theory.
Let $S_{\varOmega}$, $F_{\varOmega}$, $\Delta_{\varOmega}$
and $J_{\varOmega}$
be the modular operations of the pair $\{M,\varOmega\}$ with their usual
meanings, and which
are associated to the respective actions of $M$ and $M^{\,\prime}$,
with their respective dense standard domains of definition
${\mathcal D}\bigl(S_{\varOmega}\bigr)$,
${\mathcal D}\bigl(F_{\varOmega}\bigr)$ and ${\mathcal
D}\bigl({\Delta}_{\varOmega}\bigr)$ within ${\mathcal H}$. Remind
that $S_{\varOmega}$ and $F_{\varOmega}$ are the closures of
the (usually unbounded but closable) antilinear operators which arise from the
action of the hermitian conjugation
$M\varOmega\ni x\varOmega\,\longmapsto\,x^*\varOmega$ and
$M^{\,\prime}\varOmega\ni
y\varOmega\,\longmapsto\,y^*\varOmega$ on $M$ and
$M^{\,\prime}$, respectively. The polar decomposition of
$S_{\varOmega}$ reads as
$S_{\varOmega}=J_{\varOmega}\Delta_{\varOmega}^{1/2}$, where $\Delta_{\varOmega}=
S_{\varOmega}^*S_{\varOmega}$ is the modular operator, which is linear,
and the modular conjugation
$J_{\varOmega}$, which is antiunitary and selfadjoint, and therefore obeys
$J_{\varOmega}^2={\mathbf 1}$.

An important feature is that by $M\ni
x\,\longmapsto\,J_{\varOmega}xJ_{\varOmega}\in {
M}^{\,\prime}$ a $^*$-antiisomorphism between ${
M}$ and $M^{\,\prime}$ is given,
$J_{\varOmega}MJ_{\varOmega}=M^{\,\prime}$
(Tomita's theorem). Especially, for
each $\varphi\in {\mathcal H}$ one therefore has $p(\varphi){\mathcal H}=
\overline{J_{\varOmega}MJ_{\varOmega}\varphi}=
J_{\varOmega}\overline{M\{J_{\varOmega}\varphi\}}$. Thus, by idempotency of $J_{\varOmega}$ and Tomita's theorem,
the following intertwining relation is seen\,:
\begin{equation}\label{intertwin}
p(\varphi)J_{\varOmega}=J_{\varOmega}
p'\bigl(J_{\varOmega}\varphi\bigr).
\end{equation}
Through the following setting
a self-dual cone in ${\mathcal H}$ can be associated
with $\{ M, { {\varOmega}} \}$\,:
\begin{subequations}\label{coneprop}
\begin{equation}\label{modcone}
{\mathcal P}_{{
{\varOmega}}}^{\natural}=\overline{\bigl\{
xJ_{\varOmega}xJ_{\varOmega}\varOmega\,:\,x\in M\bigr\}}=\overline{\bigl\{
yJ_{\varOmega}yJ_{\varOmega}\varOmega\,:\,y\in M^{\,\prime}\bigr\}}\,.
\end{equation}
Recall that self-duality means that the following holds\,:
$${\mathcal P}_{{{\varOmega}}}^{\natural}=\bigl\{\psi\in {\mathcal H}:\langle\psi,\xi
\rangle\geq 0,\,\forall\,\xi
\in{\mathcal P}_{{{\varOmega}}}^{\natural}\bigr\}\,.$$ 
The elements of this natural positive cone
are (pointwise) invariant under the action of
the modular conjugation $J_{\varOmega}$, and it is known from modular
theory that ${\mathcal P}_{{{\varOmega}}}^{\natural}$ (the modular cone) is generating for the real
linear space ${\mathcal H}_{\mathrm{sa}}^\varOmega$ of all
vectors fixed under the action of the modular conjugation:
\begin{equation}\label{fix}
{\mathcal H}_{\mathrm{sa}}^\varOmega=\bigl\{\psi\in {\mathcal H}: J_\varOmega\psi=\psi\bigr\}=
{\mathcal P}_{{{\varOmega}}}^{\natural}-{\mathcal P}_{{{\varOmega}}}^{\natural}.
\end{equation}
\end{subequations}
Remark that each vector $\psi\in {\mathcal H}_{\mathrm{sa}}^\varOmega$ can be
uniquely
decomposed as $\psi=\xi_+-\xi_-$,
with $\xi_+,\xi_-\in {\mathcal P}_{{{\varOmega}}}^{\natural}$ and $\xi_+\perp\,\xi_-$.

Most importantly,
to each $\nu\in M_{*+}$ there is exactly one
vector $\xi_{\nu} \in {{\mathcal P}_{{ {\varOmega}}}^{\natural}}\,\cap\,
{\mathcal S}_M(\nu)$, and it is known since the work of
\cite{Arak:74,Conn:74,Haage:75} that the mapping
\begin{subequations}\label{7}
\begin{equation}\label{7.00}
{M_{*+}} \ni {\nu} \longmapsto {\xi_{\nu}} \in
{{\mathcal P}_{{ {\varOmega}}}^{\natural}}
\end{equation}
is onto and obeys the estimate
\begin{equation}\label{7.0}
\bigl\|\xi_\nu-\xi_\varrho\bigr\|^2\leq \|\nu-\varrho\|_1\leq
\bigl\|\xi_\nu-\xi_\varrho\bigr\|\cdot\bigl\|\xi_\nu+\xi_\varrho\bigr\|\,,
\end{equation}
\end{subequations}
for any two $\nu,\varrho\in M_{*+}$,
which shows that both cones $M_{*+}$ and
${\mathcal P}_{{ {\varOmega}}}^{\natural}$ are mutually homeomorphic when
considered with
the topologies which are induced by the respective uniform
topologies of $M_*$ and ${\mathcal H}$, respectively.

Finally, remark that in modular theory, by basic facts known from \cite{Arak:74}, 
to each two cyclic and separating vectors
$\varOmega,\varOmega'$ of $M$ there is a unique
$U(\varOmega',\varOmega)\in {\mathcal U}(M^{\,\prime})$ such that
\begin{subequations}\label{ch}
\begin{equation}\label{ch.3}
{\mathcal P}_{\varOmega'}^\natural=U(\varOmega',\varOmega){\mathcal P}_\varOmega^\natural\,.
\end{equation}
Thereby, $U(\varOmega',\varOmega)={\mathbf 1}$ is fulfilled if, and only if, $\varOmega'\in {\mathcal P}_\varOmega^\natural$ holds. and the unitary $U(\varOmega',\varOmega)$ obeys the 
following chain rule\,:
\begin{equation}\label{ch.3a}
U(w\varOmega',\varOmega)=wU(\varOmega',\varOmega),\ \forall\,w\in {\mathcal U}(M^{\,\prime})\,.
\end{equation}
\end{subequations}
\begin{example}\label{ex2}
Let ${\mathsf B}({\mathcal H})$ act by multiplication
from the left on Hilbert-Schmidt operators ${\mathcal K}$ over
the separable Hilbert space ${\mathcal H}$ (call the resulting $vN$-algebra $M$).
Let $\tau$ be the unique normal trace induced on $M$ by the standard trace
`${\mathrm{tr}}$'. Then, since ${\mathcal K}$ is already complete
under the Hilbert-Schmidt norm $\|\cdot\|_2$, which corresponds to the
scalar product $\langle x,y\rangle={\mathrm{tr}\,} y^*x$, for $x,y\in{\mathcal K}$,
in terms of Example \ref{ex1}, ${\mathcal K}$ can be identified with
$L^2(M,\tau)={\mathcal L}^2(M,\tau)$, and
the usual trace class operators uniquely correspond to $L^1(M,\tau)={\mathcal L}^1(M,\tau)$, see \cite{Scha:70}. 
$\varOmega\in {\mathcal K}$ is cyclic and separating for $M$ if, and only if,
$s(|\varOmega|)=s(|\varOmega^*|)={\mathbf 1}$. In particular, each positive
Hilbert-Schmidt operator with full support can be taken for $\varOmega$. Also, each two $\nu,\varrho\in M_{*+}$ then
have the form $\nu=\tau_a$ and
$\varrho=\tau_c$, with unique
$a,c\in L^1(M,\tau)_+$, and in sense of the above identifications
the following holds with respect to each $\varOmega>{\mathbf 0}$\,:
\begin{subequations}\label{ex20}
\begin{equation}\label{ex21}
{{\mathcal P}_{{ {\varOmega}}}^{\natural}}={\mathcal L}^2(M,\tau)_+\,,\text{ with }\xi_\mu=\sqrt{m}\text{ whenever }\mu=\tau_m,\,m\in L^1(M,\tau)_+\,;
\end{equation}
\begin{equation}\label{ex22}
\|\xi_\nu-\xi_\varrho\|^2=
\|\nu\|_1+\|\varrho\|_1 -2\tau(\sqrt{a}\sqrt{c})=\|\sqrt{a}-\sqrt{c}\,\|_2^2\,.
\end{equation}
\end{subequations}
\end{example}
\subsubsection{Bures distance and modular vectors}\label{impmod}
The most interesting information arising from specializing Theorem \ref{sesq*}
to the modular context will be
that the hypothesis of Lemma \ref{sesq} can be fulfilled, for each modular vector $\varphi_0$,
that is, whenever $\varphi_0\in {\mathcal P}_\varOmega^\natural$ holds, for some
cyclic and separating vector $\varOmega\in {\mathcal H}$ of $M$. In fact, if $\{M,\varOmega\}$ is a
standard form algebra over ${\mathcal H}$, then the following is true\,:
\begin{lemma}\label{impmod0}
$\ \forall\,\nu,\varrho\in M_{*+}:\ p^{\,\prime}(\xi_{\varrho^\perp})\leq s(|h_{\xi_\varrho,\xi_\nu}|)^\perp\,.$
\end{lemma}
\begin{proof}
Let $\varrho^\perp$ be defined as in Theorem \ref{ortho}. Then, by
$\nu\perp\varrho^\perp$ one especially has $p(\xi_{\varrho^\perp})\leq
p(\xi_\nu)^\perp$. Hence, in view of \eqref{fix} the
intertwining relations \eqref{intertwin} can be applied and show that
$J_\varOmega
p\bigl({\xi}_{\varrho^\perp}\bigr)J_\varOmega=
p'\bigl({\xi}_{\varrho^\perp}\bigr)$ and
$J_\varOmega p\bigl({\xi}_{\nu}\bigr)J_\varOmega=p'(\xi_\nu)$ are fulfilled, and therefore
in view of the above $p'(\xi_{\varrho^\perp})\leq p'(\xi_\nu)^\perp$ has to hold.
On the other hand, we have $s\bigl(\bigl|h_{\xi_\varrho,\xi_\nu}\bigr|\bigr)\leq p'(\xi_\nu)$,
by \eqref{sub1}.
The assertion now will follow upon combining together the last two relations.
\end{proof}
In line with formula \eqref{impmod.1*}, for given $\{M,\varOmega\}$ consider the vector-valued map
\begin{subequations}\label{impmodi}
\begin{equation}\label{impmodi1}
M_{*+}\times\,M_{*+}\ni \{\nu,\varrho\}\,\longmapsto\,
\psi_\varOmega^\nu(\varrho)\in {\mathcal S}_M(\varrho)\subset {\cal H}
\end{equation}
which in accordance with Theorem \ref{sesq*} can be defined by
\begin{equation}\label{impmodi2}
\psi_\varOmega^\nu(\varrho)=v_{\xi_\varrho,\xi_\nu}^*\xi_\varrho+\xi_{\varrho^\perp}=
v_{\psi,\xi_\nu}^*\psi+\xi_{\varrho^\perp},\,\forall\,\psi\in {\mathcal S}_M(\varrho)\,.
\end{equation}
\end{subequations}
In accordance with Theorem \ref{sesq*} the hypothesis of
Theorem \ref{positiv1} in the modular context then may be supplemented by the following additional information\,:
\begin{theorem}[\cite{Albe:98}]\label{mod}
$\forall\,\nu,\varrho\in M_{*+}:\ d_{\mathrm{B}}(\nu,\varrho)\ =\bigl\|\psi_\varOmega^\nu(\varrho)-\xi_\nu\bigr\|\,.$
\end{theorem}
Much about structure and meaning of the map \eqref{impmodi1} is known. In context of the following result
we only mention the probably most important properties.
\begin{lemma}\label{mod1}
For each $\nu,\varrho\in M_{*+}$ and cyclic and separating vector
$\varOmega\in {\mathcal H}$ the following
properties hold\textup{:}
\begin{enumerate}
\item\label{red.3}
$\psi_\varOmega^\nu(\varrho)$ is the unique $\tilde{\psi}\in{\mathcal S}_M(\varrho)$ with
$d_{\mathrm{B}}(\nu,\varrho)=\bigl\|\tilde{\psi}-\xi_\nu\bigr\|\ \Longleftrightarrow\ \varrho^\perp=0\,;$
\item\label{red.1}
$\psi_\varOmega^\nu(\varrho)=\psi_\varOmega^\nu(\varrho-{\varrho}^\perp)+\xi_{{\varrho}^\perp}\,;$
\item\label{red.2}
$\psi_\varOmega^\nu(\varrho)\in {\mathcal H}_{\mathrm{sa}}^\varOmega\ \Leftrightarrow\ \psi_\varOmega^\nu(\varrho)\in {\mathcal P}_\varOmega^\natural\ \Leftrightarrow\ \psi_\varOmega^\nu(\varrho-{\varrho}^\perp)\in {\mathcal P}_\varOmega^\natural\ \Leftrightarrow\ \psi_\varOmega^\nu(\varrho-{\varrho}^\perp)
\in {\mathcal H}_{\mathrm{sa}}^\varOmega\,.$
\end{enumerate}
\end{lemma}
\begin{proof}
In the following, we let
$\psi=\psi_\varOmega^\nu(\varrho)$, $h=h_{\psi,\xi_\nu}$.
To prove \eqref{red.3}, note first that, for $\tilde{\psi}$ in the fibre of $\varrho$,
with $d_{\mathrm{B}}(\nu,\varrho)=\bigl\|\tilde{\psi}-\xi_\nu\bigr\|$, by
Theorem \ref{bas5}\,\eqref{bas53}, also $\tilde{h}=h_{\tilde{\psi},\xi_\nu}\geq 0$. Also, for
$w\in M^{\,\prime}$ with $w^*w=p^{\,\prime}(\psi)$ and $\tilde{\psi}=w\psi$ in accordance with
Lemma \ref{sesq}\,\eqref{positiv-1} we see that $h=\tilde{h}$ and $s(h)=w s(h)$ have to be fulfilled.
Now, suppose $\varrho^\perp\not=0$. Then
$s(\tilde{h})^\perp\tilde{\psi}\in {\mathcal S}_M(\varrho^\perp)$ cannot vanish. Hence,
for each $t\in [0,2\pi[$, also
$\psi_t=s(\tilde{h})\tilde{\psi}+\exp{(-{\mathsf{i}}t)}\,s(\tilde{h})^\perp\tilde{\psi}\in
{\mathcal S}_M(\varrho^\perp)$, with $h_{\psi_t,\xi_\nu}=h_{s(\tilde{h})\tilde{\psi},\xi_\nu}=\tilde{h}\geq0$.
Hence, $d_{\mathrm{B}}(\nu,\varrho)=\bigl\|\psi_t-\xi_\nu\bigr\|$ for all $t\in  [0,2\pi[$,
by Theorem \ref{bas5}\,\eqref{bas53}. For $t\not=0$ obviously $\psi_t\not=\tilde{\psi}$. By Theorem \ref{mod} this especially applies to $\tilde{\psi}=\psi$, and thus for $\varrho^\perp\not=0$ nonuniqueness
follows.
On the other hand, in case of $\varrho^\perp=0$, in view of the above and
Lemma \ref{sesq}\,\eqref{positiv0*},
one has $p^{\,\prime}(\psi)=s(h)$. The above condition on $w$ in this case yields
$p^{\,\prime}(\psi)=s(h)=w$, and therefore $\tilde{\psi}=\psi$, that is, uniqueness holds.

To see \eqref{red.1}, note first that  $\psi=s(h)\psi+\xi_{\varrho^\perp}$ holds, by \eqref{impmodi2}. Let
$\psi'=\psi_\varOmega^\nu(\varrho-{\varrho}^\perp)$ and $h'=h_{s(h)\psi,\xi_\nu}$.
Obviously we then have $h'=h\geq 0$. By Lemma \ref{sesq}\,\eqref{positiv0*}
and Lemma \ref{perp1}\,\eqref{perp1.1}, $s(h)\psi\in {\mathcal S}_M(\varrho-{\varrho}^\perp)$ and $(\varrho-{\varrho}^\perp)^{\perp'}=0$ have to be fulfilled. Owing to $h'\geq 0$,
$d_{\mathrm{B}}(\nu,\varrho-\varrho^\perp)=\bigl\|s(h)\psi-\xi_\nu\bigr\|$,
by Theorem \ref{bas5}\,\eqref{bas53}.
Owing to $(\varrho-{\varrho}^\perp)^{\perp'}=0$, \eqref{red.3} can be applied to the pair $\{\nu,\varrho-\varrho^\perp\}$, and then yields $\psi'=s(h)\psi$.
Substituting the latter into $\psi=s(h)\psi+\xi_{\varrho^\perp}$ then gives \eqref{red.1}.

Note that since ${\mathcal P}_\varOmega^\natural$ is a (convex) cone, the first two of the
`$\Longleftarrow$'-implications and the third `$\Longrightarrow$'-implication
within \eqref{red.2}, as well as the implication $\psi_\varOmega^\nu(\varrho-{\varrho}^\perp)
\in {\mathcal H}_{\mathrm{sa}}^\varOmega\ \Rightarrow\ \psi_\varOmega^\nu(\varrho)\in {\mathcal H}_{\mathrm{sa}}^\varOmega$, follow from \eqref{fix} and \eqref{red.1}. Owing to this and
${\mathcal P}_\varOmega^\natural\subset {\mathcal H}_{\mathrm{sa}}^\varOmega$, to see the remaining
other implications, it suffices to prove the implication
$\psi_\varOmega^\nu(\varrho)\in {\mathcal H}_{\mathrm{sa}}^\varOmega\ \Longrightarrow\ \psi_\varOmega^\nu(\varrho-{\varrho}^\perp)\in {\mathcal P}_\varOmega^\natural$. In line
with this, assume
$\psi\in {\mathcal H}_{\mathrm{sa}}^\varOmega$. Formula \eqref{red.1} in accordance with \eqref{fix} then
tells us that $\psi'\in {\mathcal H}_{\mathrm{sa}}^\varOmega$ holds. Also, in view of
Lemma \ref{sesq}\,\eqref{positiv0*} and $(\varrho-{\varrho}^\perp)^{\perp'}=0$, one has
$p^{\,\prime}(\psi')=s(h)\leq p^{\,\prime}(\xi_\nu)$.
Now, remind the fact mentioned on
in context of \eqref{fix} and saying that in this case $\psi'=\eta_+-\eta_-$ has
to be
fulfilled with vectors $\eta_+,\eta_-\in {\mathcal P}_{\varOmega}^{\natural}$,
with
$\eta_+ \perp\eta_-$. Let $\nu_\pm\in M_{*+}$ be the positive linear forms
implemented by
$\eta_\pm$. By uniqueness of the map of \eqref{7.00} one has $\xi_{\nu_\pm}=\eta_\pm$, and
by orthogonality $\eta_+ \perp\eta_-$ and since $\|\nu_\pm\|_1=\|\eta_\pm\|^2$
is fulfilled, from estimate \eqref{7.0} the
relation $\|\nu_+ - \nu_-\|_1= \|\nu_+\|_1 + \|\nu_-\|_1$ can
be inferred. Hence, $\nu_+\perp \nu_-$, that is $p\bigl(\eta_+\bigr)\perp p\bigl(\eta_-\bigr)$
holds. With the help of \eqref{fix} and \eqref{intertwin} then also
$p'\bigl(\eta_+\bigr)\perp p'\bigl(\eta_-\bigr)$. Hence, $\varrho-\varrho^\perp=\nu_++\nu_-$. The latter
implies $p\bigl(\eta_\pm\bigr)\leq p(\psi')$. Since by assumption $J_\varOmega\psi'=\psi'$ holds, in 
reasoning with \eqref{intertwin} once more again from this in view of the above
$p'\bigl(\eta_\pm\bigr)\leq p^{\,\prime}\bigl(\psi'\bigr)=s(h)\leq p'\bigl(\xi_\nu\bigr)$ can be seen.
But then especially $p'\bigl(\eta_-\bigr)\psi'=-\eta_-$, from which
$h_{\psi',\xi_\nu}\bigl(p'\bigl(\eta_-\bigr)\bigr)=-\langle\eta_-,\xi\rangle$ has to be followed.
Now, by Theorem \ref{mod} and Theorem \ref{bas5}\,\eqref{bas53}, when applied to the pair
$\{\nu,\varrho-\varrho^\perp\}$, positivity of $h_{\psi',\xi_\nu}$ can be seen.
In line with the previously derived then especially
$-\langle\eta_-,\xi_\nu\rangle\geq 0$ has to hold, for
$\eta_-,\xi_\nu\in{\mathcal P}_{\varOmega}^{\natural}$. Since
the scalar product between vectors of the natural positive cone has to be always a nonnegative
real, $\langle\eta_-,\xi_\nu\rangle=0$ is inferred to hold. As argued above, orthogonality among
vectors of the natural positive
cone then implies orthogonality of the associated $p$- and $p'$-projections,
respectively.
Thus especially
$p'\bigl(\eta_-\bigr)\perp p'\bigl(\xi_\nu\bigr)$.
This in view of $p'\bigl(\eta_-\bigr)\leq p'\bigl(\xi_\nu\bigr)$
implies $p'\bigl(\eta_-\bigr)={\mathbf 0}$. Hence $\eta_-=0$, and
thus $\psi'\in {\mathcal P}_{\varOmega}^{\natural}$ holds, and which is \eqref{red.2}. 
\end{proof}
Without proof mention some additional facts about \eqref{impmodi2}, 
see \cite{Albe:90,Albe:92.1,Albe:97.1,Albe:98}. 
\begin{lemma}\label{rem2.1a}
Let $\varOmega$ be a  cyclic and separating vector of $M$. Let $\nu,\varrho\in M_{*+}$, and 
let $f\in M_*$ be defined by $f=\bigl\langle(\cdot)|_M\psi_\varOmega^\nu(\varrho),\xi_{\nu}\bigr\rangle$, with $\psi_\varOmega^\nu(\varrho)$ of \textup{\eqref{impmodi2}}. The 
following properties hold\textup{:}
\begin{enumerate}
\item\label{rem21}
$\psi_\varOmega^\nu(\varrho)$ is continuous at each $\{\nu,\varrho\}$ with
$\varrho^\perp=0$\,;
\item\label{rem22}
$\psi_\varOmega^\nu(\varrho)\in {\mathcal P}_{\varOmega}^{\natural}\ \Longleftrightarrow\
f=f^* \ \Longleftrightarrow\
f\geq 0\,.$
\end{enumerate}
\end{lemma}
The normal linear form $f$ defined in the premises of the previous Lemma is a quite interesting 
structure. First of all, by definition of $f$ and 
according to Theorem \ref{mod}, Theorem \ref{bas5}\,\eqref{bas53} 
and Lemma \ref{bas3} one has 
\begin{subequations}\label{rem23b}
\begin{equation}\label{rem23ba}
f({\mathbf 1})=\bigl\langle\psi_\varOmega^\nu(\varrho),\xi_{\nu}\bigr\rangle=h_{\psi_\varOmega^\nu(\varrho),\xi_{\nu}}({\mathbf 1})=\sqrt{P_M(\nu,\varrho)}\,.
\end{equation}
Secondly, since orthogonality of $\varrho^\perp$ to $\nu$ also implies $p^{\,\prime}(\xi_{\varrho^\perp})\perp  p^{\,\prime}(\xi_\nu)$, from \eqref{impmodi2} especially also follows that $f$ can be rewritten as following\,:
\begin{equation}\label{rem23bb}
f=\bigl\langle(\cdot)|_M\psi_\varOmega^\nu(\varrho),\xi_{\nu}\bigr\rangle=
\bigl\langle(\cdot)|_M v_{\xi_\varrho,\xi_\nu}^*\xi_\varrho,\xi_{\nu}\bigr\rangle=
\bigl\langle(\cdot)|_M v_{\xi_\varrho,\xi_\nu}^*\xi_\varrho,s(|h_{\xi_\nu,\xi_\varrho}^*|)
\xi_{\nu}\bigr\rangle\,,  
\end{equation}
\end{subequations}
where also \eqref{sub1} and the relation $h_{\xi_\nu,\xi_\varrho}^*=h_{\xi_\varrho,\xi_\nu}$ have 
been taken into account. Note that according to Lemma \ref{sesq}\,\eqref{positiv0*}, 
$ v_{\xi_\varrho,\xi_\nu}^*\xi_\varrho\in {\mathcal S}_M(\varrho-\varrho^\perp)$ and 
$s(|h_{\xi_\nu,\xi_\varrho}^*|)\xi_{\nu}\in {\mathcal S}_M(\nu-\nu^\perp)$ hold. 
But then, \eqref{rem23ba} and \eqref{rem23bb} prove that $f$ obviously obeys the conditions 
\eqref{a} and \eqref{b} in Theorem \ref{rem23}. Thus, by Theorem \ref{rem23}, 
$g=f$ is identified. 

On the other hand, by an analogous line of reasoning, $f^{\,\prime}
=\bigl\langle(\cdot)|_M\xi_{\varrho},\psi_\varOmega^\varrho(\nu)\bigr\rangle$ is obeying 
conditions Theorem \ref{rem23}\,\eqref{a}--\eqref{b}. Hence, by Theorem \ref{rem23} again,  
$g=f^{\,\prime}$. We may summarize this as follows. 
\begin{lemma}\label{zentral}
For $\nu,\varrho\in M_{*+}$, the unique $g\in M^*$ given by \textup{Theorem \ref{rem23}} is  
$$g=f=\bigl\langle(\cdot)|_M\psi_\varOmega^\nu(\varrho),\xi_{\nu}\bigr\rangle=\bigl\langle(\cdot)|_M\xi_{\varrho},\psi_\varOmega^\varrho(\nu)\bigr\rangle\,.$$
\end{lemma}
In particular, from this follows that $f$ for a given
pair $\{\nu,\varrho\}$
is uniquely determined by the associated
minimal pair, see Definition \ref{minpair}/Remark \ref{remminpair}.
\subsubsection{Bures distance and commutation of states}\label{impmod+}
Those pairs $\{\nu,\varrho\}$ of normal positive linear forms and fulfilling
$\psi_\varOmega^\nu(\varrho)\in {\mathcal P}_\varOmega^\natural$ deserve special interest.
\begin{definition}\label{commu}
Suppose $\nu,\varrho\in M_{*+}$, and be $\varOmega\in {\mathcal H}$ a cyclic and separating vector.
$\varrho$ is said to {\em commute} with $\nu$ if
$\psi_\varOmega^\nu(\varrho)=\xi_\varrho$ is fulfilled.
\end{definition}
We are going to show that commutation is a symmetric relation.
\begin{lemma}\label{mincommu0}
For each $\nu,\varrho\in M_{*+}$  the following are equivalent\textup{:}
\begin{enumerate}
\item\label{mincommu.10}
$\varrho$ commutes with $\nu$\,;
\item\label{mincommu.20}
$\nu$ commutes with $\varrho$\,.
\end{enumerate}
\end{lemma}
\begin{proof}
In view of the symmetry of the assertion, we may content ourselves with 
proving e.g.~that \eqref{mincommu.10} implies \eqref{mincommu.20}. Suppose \eqref{mincommu.10}. 
In line with this and 
Definition \ref{commu}, $\psi_\varOmega^\nu(\varrho)=\xi_\varrho$. By Theorem \ref{mod} and 
Theorem \ref{bas5}\,\eqref{bas53} we then have 
$h=h_{\xi_\varrho,\xi_\nu}=h_{\xi_\nu,\xi_\varrho}\geq 0$. From this with the help of 
Tomita's theorem and owing to \eqref{fix} for each $x\in M$, and corresponding $y\in M^{\,\prime}$ with $y=J_\varOmega xJ_\varOmega$, we infer that  
$f(x^*x)=\bigl\langle x^*xJ_\varOmega\psi_\varOmega^\nu(\varrho),J_\varOmega\xi_{\nu}\bigr\rangle=
\bigl\langle \xi_{\nu},J_\varOmega x^*xJ_\varOmega\psi_\varOmega^\nu(\varrho)\bigr\rangle=
\bigl\langle y^*y\xi_{\nu},\xi_\varrho\bigr\rangle=h(y^*y)\geq 0$. 
According to Lemma \ref{zentral} we then 
also have $\bigl\langle(\cdot)|_M\xi_{\varrho},\psi_\varOmega^\varrho(\nu)\bigr\rangle=f\geq 0$. 
This is the same as $\bigl\langle(\cdot)|_M\psi_\varOmega^\varrho(\nu),\xi_{\varrho}\bigr\rangle
\geq 0$. 
The conclusion of Lemma \ref{rem2.1a}\,\eqref{rem22} when applied to $\{\varrho,\nu\}$ then yields 
$\psi_\varOmega^\varrho(\nu)=\xi_\nu$, that is, 
also $\nu$ commutes with $\varrho$. 
\end{proof}
By Theorem \ref{mod}, in case of commutation
the estimate $d_{\mathrm{B}}(\nu,\varrho)\leq \|\xi_\nu-\xi_\varrho\|$ turns into an equality. The 
conjecture is that the latter behavior is equivalent to commutation
between $\nu$ and $\varrho$, generally (though this is plausible it is {\em not} an obvious matter).  
For those pairs $\{\nu,\varrho\}$ on a general $M$ however,   
which will be of interest throughout the rest of the paper, the conjecture 
can be justified easily. This will be done now      
(see also \cite[\S\,6,\,Theorem 6.15]{Albe:92.1} for the case of minimal pairs). 
\begin{lemma}\label{mincommu}
Assume $\nu,\varrho\in M_{*+}$ with either $\nu^\perp={\mathbf 0}$ or $\varrho^\perp={\mathbf 0}$. 
The following are equivalent\textup{:}
\begin{enumerate}
\item\label{mincommu.1}
$\varrho$ commutes with $\nu$\,;
\item\label{mincommu.3}
$d_{\mathrm{B}}(\nu,\varrho)=\|\xi_\nu-\xi_\varrho\|$\,.
\end{enumerate}
\end{lemma}
\begin{proof}
The implications $\eqref{mincommu.1}\,\Longrightarrow\,\eqref{mincommu.3}$ follows at once from
Definition \ref{commu} and Theorem \ref{mod}. Suppose \eqref{mincommu.3} to be fulfilled, and be 
$\varrho^\perp={\mathbf 0}$. Then, by Theorem \ref{mod} and Lemma \ref{mod1}\,\eqref{red.3}, 
$\psi_\varOmega^\nu(\varrho)=\xi_\varrho$, and \eqref{mincommu.1} is seen. In case of $\nu^\perp={\mathbf 0}$ the same reasoning by Theorem \ref{mod} and Lemma \ref{mod1}\,\eqref{red.3} for the pair 
$\{\varrho,\nu\}$ will yield $\psi_\varOmega^\varrho(\nu)=\xi_\nu$, which in view of 
Lemma \ref{mincommu0} amounts to \eqref{mincommu.1} again.
\end{proof}
On the other hand, in important special cases of $M$, like matrix algebras and, more generally, $vN$-algebras of
type ${\mathrm{I}}$, the conjecture can be shown easily, 
for each pair of normal positive linear forms.  
For the type-${\mathrm{I}}_\infty$-factor a proof will be given.
\begin{example}\label{ex3}
Suppose $M\simeq {\mathsf B}({\mathcal H})$ and $\tau$
as in Example \ref{ex2}, and suppose
$\varOmega>{\mathbf 0}$. Let $\nu=\tau_a,\,\varrho=\tau_c$, with positive
trace class operators $a$ and $c$. Suppose $\varrho$ commutes with $\nu$. By \eqref{budi},
\eqref{ex22}, \eqref{bei1} and Theorem \ref{mod}, Definition \ref{commu} then amounts to
$$\tau(\sqrt{a}\sqrt{c})=\sqrt{P_M(\nu,\varrho)}=\tau(|\sqrt{a}\sqrt{c}|)\,.$$
Since $\tau$ is faithful, by \eqref{bei2} the previous is
equivalent to  commutation $ac=ca$.

On the other hand, suppose $\nu=\tau_a,\,\varrho=\tau_c$, with positive
trace class operators $a$ and $c$ which are commuting, in the sense of operator theory.
Then, since $M^{\,\prime}$ corresponds to
right actions of bounded linear operators on Hilbert-Schmidt operators, in view of
\eqref{ex22} we have
\begin{equation*}
h(\cdot)=h_{\xi_\varrho,\xi_\nu}(\cdot)=\tau(\sqrt{a}\sqrt{c}(\cdot))\geq 0\,.
\tag{$\star$}
\end{equation*}
It follows that $s(h)$ corresponds to right multiplication with $s(\sqrt{a}\sqrt{c})=s(a) s(c)$ (remind
that for commuting $a,c$ also $s(a)$ and $s(c)$ must commute). In accordance with \eqref{deforth}
we then have that $\eta= s(h)^\perp\xi_\varrho= \sqrt{c}\, \{s(a) s(c)\}^\perp\geq {\mathbf 0}$
must be implementing for $\varrho^\perp$.
In view of \eqref{ex21} we then even have
$\xi_{\varrho^\perp}=\sqrt{c}\, \{s(a) s(c)\}^\perp$. But then, the formula \eqref{impmodi2} in view of
($\star$) reads as
$$\psi_{\varOmega}^\nu(\varrho)=s(h)\xi_\varrho +\xi_{\varrho^\perp}=\sqrt{c}\, s(a)s(c)+
\sqrt{c}\, \{s(a) s(c)\}^\perp=\sqrt{c}=\xi_\varrho\,,$$
that is, $\varrho$ commutes with $\nu$ in the sense of Definition \ref{commu}.
Hence, $ac=ca$ is equivalent to  commutation of $\varrho$ with $\nu$. Since according to
\eqref{ex22}, \eqref{bei2} and as remarked above  $ac=ca$ is also equivalent to 
$d_{\mathrm{B}}(\nu,\varrho)=\|\xi_\nu-\xi_\varrho\|$, and since by Lemma \ref{mincommu0} 
commutation is symmetric in
$\nu$ and $\varrho$, we may summarize as follows\,:

\noindent
{\em{
For $M\simeq {\mathsf B}({\mathcal H})$ and  normal positive linear forms
$\nu=\tau_a$ and $\varrho=\tau_c$, commutation of $\nu$ with $\varrho$ (resp.~$\varrho$ with $\nu$)
in the sense of
Definition \ref{commu} is equivalent to commutation of the operator densities $a$ and $c$, and
both are also equivalent to the occurence of the relation
$d_{\mathrm{B}}(\nu,\varrho)=\|\sqrt{a}-\sqrt{c}\,\|_2=\|\xi_\nu-\xi_\varrho\|$.
}}

Finally, note that owing to eqs.~\eqref{ch} in the case at
hand we must have equivalence of commutation to
$d_{\mathrm{B}}(\nu,\varrho)=\|\xi_\nu-\xi_\varrho\|$, for {\em each}
cyclic and separating $\varOmega$ (and not just for only those with
$\varOmega>{\mathbf 0}$). This especially proves that the assertion of
Lemma \ref{mincommu} remains true, for {\em each} pair $\{\nu,\varrho\}$ of
normal positive linear forms
on a factor $M$ of type ${\mathrm{I}}_\infty$ (finite factors were dealt with in
\cite{Albe:92.1}).
\end{example}
Without proof mention yet some generalities on `commutation', see 
\cite{Albe:90,Albe:92.1,Albe:97.1,Albe:98}.
\begin{remark}\label{rem2}
\begin{enumerate}
\item\label{rem2.1}
Following \cite[p.\,325]{Albe:92.1}, in a ${\mathsf C}^*$-algebraic context
`commutation' might be defined by requiring $g\geq 0$, for $g$ defined in 
Theorem \ref{rem23}. In line with Lemma \ref{rem2.1a}\,\eqref{rem22}
and Lemma \ref{zentral}, in the standard form case this then 
amounts to Definition \ref{commu}.
\item\label{rem2.2}
Note that the notion of commutation 
also does not depend on the special $\varOmega$ used (cf.~\eqref{ch.3} for the reasons of
this covariance, see \cite{Arak:74}).
It extends to arbitrary normal positive linear forms
the notion of `commutation' in the faithful case, and which refers to commutation of the respective modular
$^*$-automorphism groups, see \cite{HeTa:70,PeTa:73,Stei:93}.
\item\label{rem2.4}
According to Lemma \ref{mincommu0} and Lemma \ref{mod1}\,\eqref{red.2} one infers
that commutation must be equivalent to  commutation between the positive linear
forms of the
associated minimal pair.
Note that from Lemma \ref{mod1}\,\eqref{red.1} it follows that both
$\varrho=(\varrho-\varrho^\perp)+\varrho^\perp$ and
$\nu=(\nu-\nu^\perp)+\nu^\perp$ have to be orthogonal decompositions, in this situation.
\item\label{rem2.3}
By standard facts originating from \cite{Arak:74} it is not hard to see that, for given
$\nu,\varrho\in M_{*+}$,
$j(\nu|\varrho)=j^\varOmega(\nu|\varrho)=\frac{1}{2}
\|J_\varOmega \psi_\varOmega^\varrho(\nu)-\psi_\varOmega^\varrho(\nu)\|^2$ proves to be
independent from the special choice of the cyclic and separating vector $\varOmega$.
But then, in view of \eqref{fix}, 
Lemma \ref{mod1}\,\eqref{red.2} and Lemma \ref{mincommu0},
$j(\nu|\varrho)=0$ (resp.~$j(\varrho|\nu)=0$) if, and only if, $\nu$ commutes with $\varrho$.
Thus, $j(\nu|\varrho)$ has some invariant meaning
which in view of the
definition of the
above notion of commutation might serve as a quantitative measure estimating
how far from `commuting' with $\varrho$ the form $\nu$ is.
\end{enumerate}
\end{remark}
\subsubsection{A characterization of ${\mathcal S}_M(\nu|\varrho)$}\label{struc}
We are going to characterize ${\mathcal S}_M(\nu|\varrho)$
for a standard form $vN$-algebra in terms of either extension
properties of partial isometries of $M$, or by
comparison properties of certain characteristic orthoprojections in respect of 
the Murray-von Neumann comparability relation `$\succ$', respectively. 
\begin{theorem}\label{relfaser}
Let $\nu,\varrho\in M_{*+}$, and be $\psi=\psi_\varOmega^\nu(\varrho)$, $h=h_{\psi_\varOmega^\nu(\varrho),\xi_\nu}$. Then, $h\geq 0$, and for $u\in M^{\,\prime}$ with $u^*u=p^{\,\prime}(\xi_\nu)$
the following conditions are mutually equivalent\textup{:}
\begin{enumerate}
\item\label{relfaser.1}
$u\xi_\nu\in {\mathcal S}_M(\nu|\varrho)\,;$
\item\label{relfaser.2}
$\exists\,w\in M^{\,\prime},\,w^*w=p^{\,\prime}(\psi):\
u\,s(h)=w\,s(h)\,;$
\item\label{relfaser.3}
$us(h){u^*}^\perp \succ p^{\,\prime}(\psi)-s(h)\,.$
\end{enumerate}
\end{theorem}
\begin{proof}
By
Theorem \ref{mod} and Theorem \ref{bas5}\,\eqref{bas53} one has $h\geq 0$.
Thus $v_{\psi,\xi_\nu}=s(h)$.
By the last part of Theorem \ref{sesq*} and \eqref{impmod.1*} we then see
$\psi=s(h)\psi+\xi_{\varrho^\perp}$. Hence $s(h)^\perp\psi=\xi_{\varrho^\perp}$.
Also, since $|h_{\xi_\nu,\psi}^*|=|h_{\psi,\xi_\nu}|=|h|=h$ holds,
Lemma \ref{sesq}\,\eqref{positiv0*} may be
applied with respect to $\nu^\perp$ and then yields
$s(h)^\perp\xi_\nu\in {\mathcal S}_M(\nu^\perp)$. Note that in accordance with Theorem \ref{ortho}
especially also $\nu^\perp\perp\varrho^\perp$ must be fulfilled. This condition
is equivalent to  $\langle x s(h)^\perp\psi,ys(h)^\perp\xi_\nu\rangle=0$, for all
$x,y\in M^{\,\prime}$. Now, let $u,w\in M^{\,\prime}$ be chosen in accordance with
\eqref{relfaser.2}. By the previous $s(h)^\perp\psi\perp s(h)^\perp\xi_\nu$ and
$w s(h)^\perp\psi\perp us(h)^\perp\xi_\nu$ have to be fulfilled.
Also, owing to  $w^*w\geq s(h)$ and $u^*u\geq s(h)$ and by
\eqref{relfaser.2} the following orthogonality relations can be followed: $ws(h)\{\psi-\xi_\nu\}\perp
w s(h)^\perp\psi$ and $ws(h)\{\psi-\xi_\nu\}(=us(h)\{\psi-\xi_\nu\})\perp us(h)^\perp\xi_\nu$.
Hence, both $\{s(h)\{\psi-\xi_\nu\},s(h)^\perp\psi,s(h)^\perp\xi_\nu\}$
and $\{ws(h)\{\psi-\xi_\nu\},ws(h)^\perp\psi,us(h)^\perp\xi_\nu\}$ are orthogonal
systems of vectors, and therefore we may conclude as follows:
\begin{equation*}
\begin{split}
\|w\psi-u\xi_\nu\|^2 & =\|ws(h)\{\psi-\xi_\nu\}+ w s(h)^\perp\psi-us(h)^\perp\xi_\nu\|^2\\
& =
\|ws(h)\{\psi-\xi_\nu\}\|^2+\| w s(h)^\perp\psi\|^2+\|us(h)^\perp\xi_\nu\|^2\\
& \leq
\|s(h)\{\psi-\xi_\nu\}\|^2+\| s(h)^\perp\psi\|^2+\|s(h)^\perp\xi_\nu\|^2\\
& =
\|s(h)\{\psi-\xi_\nu\}+ s(h)^\perp\psi-s(h)^\perp\xi_\nu\|^2\\
& =\|\psi-\xi_\nu\|^2.
\end{split}
\end{equation*}
Hence, in view of Theorem \ref{mod} $\|w\psi-u\xi_\nu\|\leq d_{\mathrm{B}}(\nu,\varrho)$ must be
fulfilled. Since $w\psi\in {\mathcal S}_M(\varrho)$ and
$u\xi_\nu\in {\mathcal S}_M(\nu)$ hold, from this by Theorem \ref{bas5} equality is inferred,
$d_{\mathrm{B}}(\nu,\varrho)=\|w\psi-u\xi_\nu\|$. That is, \eqref{relfaser.1} has to be true.

Suppose \eqref{relfaser.1} is fulfilled. By definition of ${\mathcal S}_M(\nu|\varrho)$ and in view of
\eqref{faser} there
has to exist $w\in M^{\,\prime}$ with $w^*w=p^{\,\prime}(\psi)$ such that $d_{\mathrm{B}}(\nu,\varrho)=\|w\psi-u\xi_\nu\|$ holds. In view of Theorem \ref{bas5}\,\eqref{bas53} then $h_{w\psi,u\xi_\nu}\geq 0$.
Hence, $h(u^*(\cdot)w)\geq 0$, on $ M^{\,\prime}$. From $h\geq 0$ with
$s(h)\leq p^{\,\prime}(\xi_\nu)$ we then infer that $h_{w\psi,u\xi_\nu}=
h(u^*(\cdot)w)=h(u^*(\cdot)wu^*u)\geq 0$. Let $q\in  M^{\,\prime}$ be the orthoprojection
$q=u s(h) u^*$. Then $qu=u s(h)$, and thus the previous formula reads as $h_{w\psi,u\xi_\nu}
=h(u^*(\cdot)(wu^*q)u)=g((\cdot)wu^*q)\geq 0$, with $g=h(u^*(\cdot)u)\geq 0$. Note that
$s(g)=q$, and $quw^*wu^*q=u s(h)p^{\,\prime}(\psi) s(h)u^*=u s(h) u^*=q$. Thus $wu^*q$ is
a partial isometry in $ M^{\,\prime}$ with initial projection $s(g)$ and obeying
$g((\cdot)wu^*q)\geq 0$.
By uniqueness of the polar decomposition then $wu^*q=s(g)=q$ follows. By definition of $q$ from
this $w s(h)= w u^* q u=q u=u s(h)$ follows, which is the condition in \eqref{relfaser.2}.
Thus, in view of the above \eqref{relfaser.1} is equivalent to  \eqref{relfaser.2}.

Note that \eqref{relfaser.3} is equivalent to  $v^*v=p^{\,\prime}(\psi)-s(h)$ and
$vv^*\leq us(h){u^*}^\perp$, for some $v\in M^{\,\prime}$. This is the same as requiring
$v^*v=p^{\,\prime}(\psi)-s(h)$ with $v^*us(h)={\mathbf 0}$.
In defining $w=us(h)+v$ we then have $w^*w=s(h)u^*us(h)+
p^{\,\prime}(\psi)-s(h)=s(h)+p^{\,\prime}(\psi)-s(h)=p^{\,\prime}(\psi)$, with
$ws(h)=us(h)$. On the other hand, for each $w$ with $w^*w=p^{\,\prime}(\psi)$ and
$ws(h)=us(h)$ one has that $w\{p^{\,\prime}(\psi)-s(h)\}w^*\perp us(h)u^*$. Hence,
$v=w\{p^{\,\prime}(\psi)-s(h)\}$ is a partial isometry, with $v^*v=p^{\,\prime}(\psi)-s(h)$ and
$vv^*\leq us(h){u^*}^\perp$. Thus \eqref{relfaser.3} is equivalent to  \eqref{relfaser.2}.
\end{proof}
\begin{remark}\label{rem40}
Suppose $M$ is a unital ${\mathsf C}^*$-algebra, and $\pi$ is a
unital $^*$-representation such that the $\pi$-fibres of given two positive linear forms
$\nu,\varrho$ exist.
Then, due to the general character of Theorem \ref{sesq*} and Lemma \ref{sesq}, and since
owing to Theorem \ref{positiv1} in the $\pi$-fibre of $\nu$ a vector $\varphi_0$ and corresponding
$\eta$ and associated $\psi_0$ obeying eqs.~\eqref{crit} exist,
Theorem \ref{relfaser} can be seen as a special case of a (less specific) result over
unital ${\mathsf C}^*$-algebras and
characterizing ${\mathcal S}_{\pi,M}(\nu|\varrho)$ in terms of either extension
properties of partial isometries of $\pi(M)^{\,\prime}$, or by
comparison properties of certain orthoprojections within the
$vN$-algebra $\pi(M)^{\,\prime\prime}$,
respectively. It is obvious from the previous proof how this more academic 
variant should read then (with $\varphi_0,\eta$ and $\psi_0$ 
instead of $\xi_\nu,\xi_{{\varrho}^\perp}$ and 
$\psi$). 
\end{remark}
\subsubsection{Modular vectors and the structure of ${\mathcal S}_M(\nu|\varrho)$}\label{struc1}
Start with some generalities about ${\mathcal S}_M(\nu)$ and ${\mathcal S}_M(\nu|\varrho)$, and 
which have to be fulfilled, for each $\nu\in M_{*+}$ and each pair $\{\nu,\varrho\}\subset M_{*+}$, 
respectively. 

Note that according to eqs.~\eqref{ch} and with respect to 
a (particular) given standard form action $\{M,\varOmega\}$ the following has to be fulfilled\,:
\begin{equation}\label{modvec}
{\mathcal U}(M^{\,\prime})\xi_\nu={\mathcal S}_M(\nu)\cap\bigl\{\cup_{\varOmega'} {\mathcal P}_{\varOmega'}^\natural\bigr\}\,.
\end{equation}
Thereby, $\varOmega'$ may extend 
over the cyclic and separating vectors of $M$ within ${\mathcal H}$. Due to \eqref{modvec} 
we will refer to ${\mathcal U}(M^{\,\prime})\xi_\nu$ as set of 
{\em modular vectors} of ${\mathcal S}_M(\nu)$. First of all, what will be done is to 
relate ${\mathcal S}_M(\nu)$ and ${\mathcal S}_M(\nu|\varrho)$ topologically 
to the subset of modular vectors of the fibre of $\nu$.   

It is plain to see that Theorem \ref{relfaser}\,\eqref{relfaser.2} 
can be satisfied by each partial isometry 
$u\in M^{\,\prime}$ which can be obtained from an isometry $v\in M^{\,\prime}$ as  
$u=v p'(\xi_\nu)$. In fact, owing to $v^*v={\mathbf 1}$ and 
$s(h)\leq  p'(\psi)\wedge p'(\xi_\nu)$, $w=v p'(\psi)$ can be taken 
there. In view of 
Theorem \ref{relfaser}, and since unitaries are special    
isometries, we thus have the following net of inclusions to hold, 
in each case of $\nu,\varrho\in M_{*+}$\,:
\begin{equation}\label{inkl.1}
{\mathcal U}(M^{\,\prime})\xi_\nu\subset \bigl\{w\xi_\nu:\,w^*w={\mathbf 1},\,w\in M^{\,\prime}\bigr\}\subset {\mathcal S}_M(\nu|\varrho)\subset {\mathcal S}_M(\nu)\,.
\end{equation}
In case of finite $M^{\,\prime}$ (which occurs iff $M$ is finite) 
equality has to occur throughout, 
since then the whole fibre 
is made of modular vectors, see Remark \ref{fin}. On the other hand, in case of properly infinite 
$M^{\,\prime}$ by standard 
facts it is known that the strong (operator) closure of the 
unitary group ${\mathcal U}(M^{\,\prime})$ is the 
set of isometries of $M^{\,\prime}$. But then, since a general nonfinite  
$M^{\,\prime}$ can be (centrally) decomposed into two parts which are finite and 
properly infinite, respectively, the above facts relating the finite and properly infinite 
case of $M^{\,\prime}$ and \eqref{inkl.1} fit together into the following result 
($\overline{\phantom{a}}^{\,\mathrm{st}}$ indicates the closure 
with respect to the strong operator topology on $M$)\,:
\begin{subequations}\label{inkl}
\begin{equation}\label{inkl.2}
{\mathcal U}(M^{\,\prime})\xi_\nu\subset \bigl\{w\xi_\nu:\,w^*w={\mathbf 1},\,w\in M^{\,\prime}\bigr\} =\overline{{\mathcal U}(M^{\,\prime})}^{\,\mathrm{st}}\xi_\nu\subset {\mathcal S}_M(\nu|\varrho)\,.
\end{equation}
On the other hand, from \eqref{inkl.1} and by closedness of the fibre we also infer  
\begin{equation}\label{inkl.3}
\overline{{\mathcal U}(M^{\,\prime})}^{\,\mathrm{st}}\xi_\nu\subset 
\overline{{\mathcal U}(M^{\,\prime})\xi_\nu}\subset{\mathcal S}_M(\nu)\,.
\end{equation} 
Finally, for completeness, note that also the fibre itself is related to the modular vectors. 
Let, for a subset of vectors $\Gamma\subset {\mathcal H}$ and real 
$r>0$, $\Gamma[r]$ be defined by $\Gamma[r]=\{\varphi\in \Gamma:\|\varphi\|=r\}$. 
Then, the formula for ${\mathcal S}_M(\nu)$ reads as follows (`${\mathop{\mathrm{conv}}}$' stands for the operation of taking the convex hull)\,:
\begin{equation}\label{inkl.4}
{\mathcal S}_M(\nu)=\overline{{\mathop{\mathrm{conv}}}\, 
{\mathcal U}(M^{\,\prime})\xi_\nu}\left[\sqrt{\|\nu\|_1}\right]\,.
\end{equation}
\end{subequations}
In fact, according to \cite{DyRu:66}, especially  
$\bigl(M^{\,\prime}\bigr)_1=\overline{{\mathop{\mathrm{conv}}}\, {\mathcal U}(M^{\,\prime})}$ (uniform closure) must be fulfilled. Hence, in view of \eqref{bas4},  
$
{\mathcal S}_M(\nu)\subset \overline{{\mathop{\mathrm{conv}}}\, {\mathcal U}(M^{\,\prime})}\,\xi_\nu\subset \overline{{\mathop{\mathrm{conv}}}\, {\mathcal U}(M^{\,\prime})\xi_\nu}
$
can be followed. On the other hand, suppose $\varphi\in {\mathcal H}$, with 
$\|\varphi\|^2=\|\nu\|_1=\nu({\mathbf 1})$ and $\varphi=\lim_{n\to\infty} \varphi_n$, 
with a sequence $\{\varphi_n\}\subset {\mathop{\mathrm{conv}}}\, 
{\mathcal U}(M^{\,\prime})\xi_\nu$. Let us define $\nu_n\in M_{*+}$ by the condition that 
$\varphi_n\in {\mathcal S}_M(\nu_n)$. Then, by the previous especially $\varphi_n=a_n\xi_\nu$, 
with $a_n\in \bigl(M^{\,\prime}\bigr)_1$, and thus $\nu_n\leq \nu$, for each subscript. But then 
$\|\nu-\nu_n\|_1=\nu({\mathbf 1})-\nu_n({\mathbf 1})=\|\varphi\|^2-\|\varphi_n\|^2$. In view of 
the above $\nu=\lim_{n\to\infty} \nu_n$ follows, in the uniform sense. 
Hence,  
$
\overline{{\mathop{\mathrm{conv}}}\, {\mathcal U}(M^{\,\prime})\xi_\nu}[\sqrt{\|\nu\|_1}]\subset {\mathcal S}_M(\nu)
$
holds. Taking together the previous inclusions yields \eqref{inkl.4}.

Note that if $\nu$ is faithful, that is, ${\mathbf 1}=s(\nu)=p(\xi_\nu)$ is fulfilled, then owing to 
\eqref{intertwin} one also has $p^{\,\prime}(\xi_\nu)={\mathbf 1}$, and then in view of 
\eqref{bas4}, \eqref{inkl.1}   
and eqs.~\eqref{inkl} we infer  
\begin{equation}\label{inklfolg}
{\mathcal S}_M(\nu|\varrho)=\overline{{\mathcal U}(M^{\,\prime})}^{\,\mathrm{st}}\xi_\nu=\overline{{\mathcal U}(M^{\,\prime})\xi_\nu}=\overline{{\mathop{\mathrm{conv}}}\, {\mathcal U}(M^{\,\prime})\xi_\nu}\left[\sqrt{\|\nu\|_1}\right]={\mathcal S}_M(\nu)\,,
\end{equation}
for each $\varrho\in M_{*+}$. Thus, from point of view of Problem \ref{prob}\,\eqref{prob4} those 
pairs with faithful $\nu$ will not be of any interest. 

In contrast to the previous, the following observation likely yields the most important class of pairs $\{\nu,\varrho\}$ where
an easy to handle formula for ${\mathcal S}_M(\nu|\varrho)$ can be derived, and which is nontrivial 
insofar as \eqref{inklfolg} then need not be true.  

\noindent
To explain the formula, suppose in the notations of 
Theorem \ref{relfaser}, with  
$\psi=\psi_\varOmega^\nu(\varrho)$ and $h=h_{\psi_\varOmega^\nu(\varrho),\xi_\nu}$, for given 
$\{\nu,\varrho\}$ the conditions 
\begin{equation*}
p^{\,\prime}(\psi)={\mathbf 1},\ s(h)=p^{\,\prime}(\xi_\nu)
\tag{$\circ$}
\end{equation*} 
to be fulfilled. Then, by 
Theorem \ref{relfaser}\,\eqref{relfaser.1}--\eqref{relfaser.2}, for given $u\in M^{\,\prime}$ 
with $u^*u=p^{\,\prime}(\xi_\nu)$ and $u\xi_\nu\in {\mathcal S}_M(\nu|\varrho)$ we find 
$w \in M^{\,\prime}$ with $w^*w={\mathbf 1}$ such that $u=u p^{\,\prime}(\xi_\nu)=u s(h)=w s(h)=w
p^{\,\prime}(\xi_\nu)$. Hence, $u\xi_\nu=w p^{\,\prime}(\xi_\nu)\xi_\nu=w\xi_\nu$. 
Thus, in this case 
${\mathcal S}_M(\nu|\varrho)\subset \{w\xi_\nu:\,w^*w={\mathbf 1},\,w\in M^{\,\prime}\}$, and then 
in view of \eqref{inkl.2} the formula 
\begin{subequations}\label{allg0}
\begin{equation}\label{allg}
{\mathcal S}_M(\nu|\varrho)=\bigl\{w\xi_\nu:\,w^*w={\mathbf 1},\,w\in M^{\,\prime}\bigr\}=\overline{{\mathcal U}(M^{\,\prime})}^{\,\mathrm{st}}\xi_\nu
\end{equation}
can be followed. In view of \eqref{inkl.3} the closure of the latter set then reads as  
\begin{equation}\label{allg1}
\overline{{\mathcal S}_M(\nu|\varrho)}=\overline{{\mathcal U}(M^{\,\prime})\xi_\nu}\,.
\end{equation}
\end{subequations}
\begin{example}\label{main1}
The easiest way to cope with the premises of 
($\circ$) is to assume that $\varrho$ is faithful, and that $\nu$ is commuting with $\varrho$.

In fact, for faithful $\varrho$, 
by Theorem \ref{ortho}
$\nu^\perp=0$ with respect to $\varrho$, for each $\nu\in M_{*+}$. Remind that $h\geq 0$, by 
Theorem \ref{relfaser}. 
Hence, $\{p^{\,\prime}(\xi_\nu)-s(h)\}\xi_\nu=s(h)^\perp\xi_\nu=0$, 
by Lemma \ref{sesq}\,\eqref{positiv0*}. Thus
$p^{\,\prime}(\xi_\nu)=s(h)$. By assumption and Lemma \ref{mincommu0}, 
$\psi=\xi_\varrho$. Hence, in view of \eqref{intertwin} and \eqref{fix},
faithfulness of $\varrho$, which means that $\psi$ is separating, implies the modular vector 
$\psi$ to be also cyclic, $p^{\,\prime}(\psi)={\mathbf 1}$. Thus, the conditions  
of ($\circ$) are satisfied and, as mentioned above, both formulae of eqs.~\eqref{allg0} then 
must be fulfilled. 
\end{example}
\begin{remark}\label{rem4}
It is easy to see that \eqref{allg} is a special case of the formula  
\begin{equation}\label{allg.0}
{\mathcal S}_M(\nu|\varrho)=\bigl\{w\xi_\nu:\,w^*w=p^{\,\prime}(\psi),\,w\in M^{\,\prime}\bigr\}
\end{equation}
which holds provided $\nu^\perp=0$ is fulfilled. But note that in general it is difficult 
to say something definite about $p^{\,\prime}(\psi)$ without having some additional informations 
on the pair $\{\nu,\varrho\}$ and which go beyond $\nu^\perp=0$. Even for a 
faithful $\varrho$, 
which implies $\nu^\perp=0$ e.g., in general only $p^{\,\prime}(\psi)\sim {\mathbf 1}$ 
can be followed. Note that from \eqref{allg} especially follows that under the premises of 
Example \ref{main1} the modular vectors in the fibre of $\nu$ are dense within ${\mathcal S}_M(\nu|\varrho)$.    
\end{remark}
\subsubsection{When does ${\mathcal S}_M(\nu|\varrho)$ provide the whole fibre?}
We continue with conditions under which Problem \ref{prob}\,\eqref{prob3} can be
affirmatively answered, for a pair $\nu,\varrho\in M_{*+}$, and which will lead us
far beyond the finite case or the situation known from \eqref{inklfolg}.
\begin{lemma}\label{einfach0}
Let $\varrho^\perp=0$, or be $s(\nu)$ finite.
Then ${\mathcal S}_M(\nu|\varrho)={\mathcal S}_M(\nu)$.
\end{lemma}
\begin{proof}
Suppose  $\varrho^\perp=0$ first.
By Lemma \ref{sesq}\,\eqref{positiv0*} then $\{p^{\,\prime}(\psi)-s(h)\}\psi=
s(h)^\perp\psi=0$. Hence $p^{\,\prime}(\psi)=s(h)$, and thus the condition of
Theorem \ref{relfaser}\,\eqref{relfaser.3} can be fulfilled in a trivial way, for each $u\in M^{\,\prime}$,
$u^*u=p^{\,\prime}(\xi_\nu)$. In view of \eqref{bas4} the fact then follows.

If the support orthoprojection $s(\nu)=p(\xi_\nu)$ is finite, then an
application of \eqref{intertwin} yields that also
$p^{\,\prime}(\xi_\nu)$ is finite (and vice versa), and so has to be also each other
$p^{\,\prime}(\chi)$, for $\chi\in {\mathcal S}_M(\nu)$. Let
$\chi=u\xi_\nu$, $u^*u=p^{\,\prime}(\xi_\nu)$, with $u\in M^{\,\prime}$.
Since also $z=p^{\,\prime}(\xi_\nu)\vee
p^{\,\prime}(\chi)$ is finite, by \cite[2.4.2.]{Saka:71} $u\in z M^{\,\prime}z$ extends to a 
unitary $u_1\in z M^{\,\prime}z$. The partial isometry $w=(u_1+z^\perp)p^{\,\prime}(\psi)$
then obeys $u s(h)=w s(h)$ and $w^*w=p^{\,\prime}(\psi)$, and thus Theorem \ref{relfaser} can be applied.
\end{proof}
\begin{remark}\label{rem4a}
Under the same premises on $\pi$ as in Remark \ref{rem4}, there must be a ${\mathsf C}^*$-variant of
that part of Lemma \ref{einfach0} which refers to the assumption $\varrho^\perp=0$, accordingly.
This together with \eqref{perp} will show that for each such $\pi$
by the condition
${\varrho}^\perp=0$ then ${\mathcal S}_{\pi,M}(\nu|\varrho)={\mathcal S}_{\pi,M}(\nu)$ will be
implied.
\end{remark}
Another observation provides conditions under which existence of cyclic elements
in a fibre allows us to decide on whether or not ${\mathcal S}_M(\nu|\varrho)={\mathcal S}_M(\nu)$ would
hold.
\begin{lemma}\label{einfach}
Suppose $\nu,\varrho\in M_{*+}$. The following assertions are valid.
\begin{enumerate}
\item\label{einfach.1}
Suppose $\varrho^\perp\not=0$ and $\nu^\perp=0$. If a cyclic vector $\chi$
in the fibre of $\nu$ exists, then each such vector obeys $\chi\not\in
{\mathcal S}_M(\nu|\varrho)$. Hence then
${\mathcal S}_M(\nu|\varrho)\not={\mathcal S}_M(\nu)$.
\item\label{einfach.2}
In the factor case of $M$, if ${\mathcal S}_M(\nu|\varrho)\not={\mathcal S}_M(\nu)$ holds,
then $\varrho^\perp\not=0$ and there exists
a cyclic vector $\chi$ in the fibre of $\nu$.
\end{enumerate}
\end{lemma}
\begin{proof}
According to Lemma \ref{sesq}\,\eqref{positiv0*} and by positivity of $h$
the assumptions on
$\varrho^\perp$ and $\nu^\perp$ equivalently read as $\{p^{\,\prime}(\psi)-s(h)\}\psi=
s(h)^\perp\psi\not=0$ and
$\{p^{\,\prime}(\xi_\nu)-s(h)\}\xi_\nu=s(h)^\perp\xi_\nu=0$. Thus $p^{\,\prime}(\psi)>s(h)$ and
$p^{\,\prime}(\xi_\nu)=s(h)$. Suppose a cyclic $\chi$ in the fibre of $\nu$ to exist. Then, in view of
\eqref{bas4} there is $u\in M^{\,\prime}$ with $u^*u=p^{\,\prime}(\xi_\nu)$ and $\chi=u\xi_\nu$.
Thus especially $uu^*=p^{\,\prime}(\chi)={\mathbf 1}$ and $p^{\,\prime}(\psi)>s(h)$. It is obvious that
Theorem \ref{relfaser}\,\eqref{relfaser.3} cannot be satisfied by $u$.
Thus $\chi=u\xi_\nu\not\in {\mathcal S}_M(\nu|\varrho)$, which is \eqref{einfach.1}.

To see \eqref{einfach.2}, suppose now $M$ to be a factor,
and ${\mathcal S}_M(\nu|\varrho)\not={\mathcal S}_M(\nu)$.
By Lemma \ref{einfach0} we then have $\varrho^\perp\not=0$ and infinite $s(\nu)$.
Remind that, $M$ being a standard form $vN$-algebra, $M$
admits faithful normal positive linear forms. Hence $M$ is a countably decomposable factor.
Thus there can exist only one equivalence class of infinite orthoprojections in $M$, and therefore
$p(\xi_\nu)=s(\nu)\sim {\mathbf 1}$ must hold. From the latter in view of
\eqref{intertwin} and \eqref{fix} the relation
$p^{\,\prime}(\xi_\nu)\sim {\mathbf 1}$ follows. Thus, a cyclic vector has to
exist in ${\mathcal S}_M(\nu)$, and $\varrho^\perp\not=0$ by the above.
\end{proof}
Since in the factor case and with faithful $\varrho$ the condition $\nu^\perp=0$ is fulfilled for any
$\nu\in M_{*+}$, in this case Lemma \ref{einfach}\,\eqref{einfach.1} and \eqref{einfach.2} then provide
sufficient
and necessary conditions for ${\mathcal S}_M(\nu|\varrho)\not={\mathcal S}_M(\nu)$ to occur,
or equivalently, for Problem \ref{prob}\,\eqref{prob3} to have an affirmative answer. Thus, in
this case also the opposite direction of the implication in Lemma \ref{einfach0} holds true.
\begin{theorem}\label{fac}
Let $M$ be a factor which acts in standard form. For $\nu,\varrho\in M_{*+}$ with faithful
$\varrho$,
${\mathcal S}_M(\nu|\varrho)={\mathcal S}_M(\nu)$ holds if, and only if, $\varrho^\perp=0$ or $s(\nu)$ is
finite.
\end{theorem}
\begin{proof}
As mentioned, under the given premises $\nu^\perp=0$. Thus
Lemma \ref{einfach} equivalently says that ${\mathcal S}_M(\nu|\varrho)={\mathcal S}_M(\nu)$ holds iff
$\varrho^\perp=0$ or $p^{\,\prime}(\xi_\nu)\not\sim {\mathbf 1}$. The latter is equivalent to 
$s(\nu)=p(\xi_\nu)\not\sim {\mathbf 1}$, by \eqref{intertwin} and \eqref{fix}. Since $M$ is
a countably decomposable factor, this is the same as asserting $s(\nu)$ to be finite.
\end{proof}
\begin{remark}\label{rem42}
For finite dimensional or commutative $M$ a test on
$\varrho^\perp\not=0$ can be achieved easily. In fact, in all these cases
$\varrho^\perp\not=0$ proves equivalent to  the condition
$s(\varrho)\wedge s(\nu)^\perp\not={\mathbf 0}$. Thus, in these
cases for the question on $\varrho^\perp\not=0$ only the mutual
position of the respective support orthoprojections within the projection lattice
of the algebra matters, refer to \cite{Hein:91} for that.
\end{remark}
Unfortunately, in the infinite dimensional, noncommutative cases of $M$, the
condition $s(\varrho)\wedge s(\nu)^\perp\not={\mathbf 0}$ in general does not
suffice any longer to imply $\varrho^\perp\not=0$.
For convenience of the reader we include a counterexample. Thereby, in view of
\eqref{perp} for the effect to appear it cannot be of relevance whether $M$ acts in standard form or not,
and thus for simplicity e.g.~the algebra
$M={\mathsf B}({\mathcal H})$ of all bounded linear operators over some separable,
infinite-dimensional Hilbert space will be considered.
\begin{example}\label{ex}
Let $\{\varphi_k:k\in {\mathbb N}\}$ be a complete orthonormal system of vectors within ${\mathcal H}$, and
$\psi\in {\mathcal H}$ be a unit vector with $\langle \psi,\varphi_k\rangle\not=0$, for all $k\in {\mathbb N}$. Let pure normal states $\nu_k$ be defined as $\varphi_k\in {\mathcal S}_M(\nu_k)$, and for $0<\beta<1$
define $\varrho\in M_{*+}$ by
$\varrho=\sum_{k\in {\mathbb N}} \beta^k |\langle \psi,\varphi_k\rangle|^2 \nu_k$. Let
$\nu\in M_{*+}$ be any normal positive linear form with support orthoprojection $s(\nu)=p_\psi^\perp$, where
$p_\psi$ is the orthoprojection onto ${\mathbb C}\,\psi$. Thus, $\varrho$ is faithful and $\nu$ has support with codimension one. We now consider $\varrho^\perp$ in the sense of
Theorem \ref{ortho} for $\{\nu,\varrho\}$. By orthogonality with $\nu$ then
$\varrho^\perp=\gamma\cdot \nu_\psi$, where $\nu_\psi$ is the vector state of $\psi$, and
$\gamma\geq 0$, real. On the other hand, $\varrho^\perp\leq \varrho$ has to be fulfilled.
Hence, for each $k\in {\mathbb N}$ especially,
$$\gamma\,|\langle \psi,\varphi_k\rangle|^2=\varrho^\perp(p_{\varphi_k})\leq \varrho(p_{\varphi_k})
= \beta^k |\langle \psi,\varphi_k\rangle|^2 \,.$$
From this $\gamma\leq \beta^k$ follows, for all $k\in {\mathbb N}$. Thus $\gamma=0$, and $\varrho^\perp=0$
is seen. But note that $s(\varrho)\wedge s(\nu)^\perp=p_\psi$ in this case.
\end{example}

\subsubsection{Fibres and centralizers}\label{impmoda}
In line with Problem \ref{prob}\,\eqref{prob4} now we are going to look for conditions
under which on a standard
form $vN$-algebra $M$ and with fixing
some individual vector $\varphi_0$ in the fibre of $\nu$ the hypothesis of
Theorem \ref{positiv1} could be violated.
To that aim, we start looking
for additional conditions on the support of $\nu$ under which $\varrho^\perp$ can be constructed
more explicitly, in each case of $\varrho$. In making reference to the notion of the {\em centralizer
$vN$-algebra} $M^\varrho$, which reads as
$$
M^\varrho=\{x\in M:\,\varrho(xy)=\varrho(yx),\,\forall\,y\in M\}\,,$$
such conditions can be easily formulated.
\begin{lemma}\label{conperp}
Suppose $\nu,\varrho\in M_{*+}$, with $s(\nu)\in M^\varrho$.
Then $\varrho^\perp=\varrho(s(\nu)^\perp(\cdot))$.
\end{lemma}
\begin{proof}
Owing to $s(\nu)\in M^\varrho$ both $\varrho(s(\nu)^\perp(\cdot))$ and
$\varrho(s(\nu)(\cdot))$ have to be
normal positive linear forms which sum up to $\varrho$. Especially,
$\omega=\varrho(s(\nu)^\perp(\cdot))=\varrho(s(\nu)^\perp(\cdot)s(\nu)^\perp)$ then is subordinate to
$\varrho$ and is
orthogonal to $\nu$ by construction. On the other hand, $s(\varrho^\perp)\leq s(\nu)^\perp$
also implies $\varrho^\perp=\varrho^\perp(s(\nu)^\perp(\cdot)s(\nu)^\perp)$. Hence,
$\omega-\varrho^\perp=\{\varrho-\varrho^\perp\}(s(\nu)^\perp(\cdot)s(\nu)^\perp)\geq 0$.
In view of Theorem \ref{ortho} then $\omega=\varrho^\perp$ follows.
\end{proof}
Now, constructions where
Lemma \ref{einfach} can be applied, can be easily achieved.
\begin{lemma}\label{cent}
Suppose $\nu,\varrho\in M_{*+}$, with $s(\nu)\in M^\varrho$, $s(\nu)<s(\varrho)$ and
$s(\nu)\sim {\mathbf 1}$.
Then ${\mathcal S}_M(\nu|\varrho)\not={\mathcal S}_M(\nu)$.
\end{lemma}
\begin{proof}
In view of Lemma \ref{conperp} and $s(\nu)<s(\varrho)$, $\varrho^\perp=
\varrho(s(\nu)^\perp(\cdot))\not=0$ follows. On the other hand, $s(\nu)\leq s(\varrho)$
implies $\nu^\perp=0$. Moreover, in view of
\eqref{intertwin} and \eqref{fix} from $p(\xi_\nu)=s(\nu)\sim {\mathbf 1}$ the relation
$p^{\,\prime}(\xi_\nu)\sim {\mathbf 1}$ is obtained. Thus a cyclic vector has to
exist in ${\mathcal S}_M(\nu)$. Application of Lemma \ref{einfach}\,\eqref{einfach.1}
then yields the result.
\end{proof}
Clearly, in order that examples in accordance with these premises could exist
at all, $M$ has to be supposed to be infinite.
Moreover, for factors and faithful $\varrho$, the assertion of
Lemma \ref{cent} even can be strengthened.
\begin{lemma}\label{centfac}
On an infinite factor $M$, for $\nu,\varrho\in M_{*+}$ with faithful
$\varrho$, and with $s(\nu)\in M^\varrho$,
${\mathcal S}_M(\nu|\varrho)\not={\mathcal S}_M(\nu)$ happens if, and only if, $s(\nu)$ is
infinite and $s(\nu)\not={\mathbf 1}$.
\end{lemma}
\begin{proof}
As in the previous proof, $s(\nu)\not={\mathbf 1}$ and $s(\nu)\in M^\varrho$ for faithful $\varrho$
imply $\varrho^\perp\not=0$. Note that the latter condition implies $s(\nu)\not={\mathbf 1}$, by
triviality. Hence, for faithful $\varrho$ and $s(\nu)\in M^\varrho$, $s(\nu)\not={\mathbf 1}$ is
equivalent to  $\varrho^\perp\not=0$. But then the assertion
is a consequence of Theorem \ref{fac}.
\end{proof}
\noindent
In order to see along the lines of Lemma \ref{cent}/\ref{centfac} that 
${\mathcal S}_M(\nu|\varrho)\not={\mathcal S}_M(\nu)$ can happen, 
first recall some further facts from modular and operator theory.
\begin{remark}\label{rem3}
\begin{enumerate}
\item\label{rem32}
In case of a standard form $vN$-algebra $\{M,\varOmega\}$, $\varrho$ defined by
$\varOmega\in {\mathcal S}_M(\varrho)$ is faithful. A fundamental result of
modular theory then says that $M^\varrho$ is the same as the
fixed-point algebra of the modular $^*$-automorphism group of $\varrho$\,:
$$
{M}^\varrho=\bigl\{x\in {M}:{\varSigma}_t^\varrho(x)=x,\ \forall\, t\in{\mathbb{R}}\bigr\}\,,
$$
where the {\em modular $^*$-automorphism group  $\{{\varSigma}_t^\varrho\}$ of $\varrho$} for all $t\in {\mathbb{R}}$
acts on elements $x\in M$ as ${\varSigma}_t^\varrho(x)=
\Delta_{\varOmega}^{{\mathsf i}\,t}x
\Delta_{\varOmega}^{-{\mathsf i}\,t}$. As consequence of this
one has
\begin{subequations}\label{add}
\begin{equation}\label{add.2}
\overline{\bigl(M^\varrho\bigr)_{\mathrm{h}}\varOmega}=
\overline{M_{\mathrm{h}}\varOmega}\cap
\overline{M^{\,\prime}_{\mathrm{h}}\varOmega}\,,
\end{equation}
\begin{equation}\label{add.1}
\overline{\bigl(M^\varrho\bigr)_{\mathrm{h}}\varOmega}\subset {\mathcal H}_{\mathrm{sa}}^\varOmega\,,
\end{equation}
with ${\mathcal H}_{\mathrm{sa}}^\varOmega$ of \eqref{fix}, and where $(\cdot)_{\mathrm{h}}$
refers to the Hermitian portion of the algebra in question. For these facts we refer the reader
to \cite{KoSz:66,Take:72}. Remark that \eqref{add.2} remains true with positive portions instead
of Hermitian parts, accordingly. \end{subequations}
\item\label{rem33}
Relating the structure of $M^\varrho$, if $M$ is properly infinite,
by the duality theorem of Connes and
Takesaki \cite{Conn:73,Take:73}, see also \cite{Naka:77}, there are a semifinite but infinite-dimensional 
$vN$-algebra $N$ and a $\sigma$-weakly continuous group $\bigl(\alpha_t\bigr)_{t\in{\mathbb{R}}}$
of $^*$-automorphisms of $N$
such that ${M}$ is isomorphic to a so-called
${\mathsf W}^*$-crossed product, $ N\,\bigotimes_\alpha{\mathbb{R}}\simeq{
M}$. Thereby, one knows that $N$ can be
isomorphically identified with the fixed-point algebra of $\{{\varSigma}_t^\varrho\}$.
Hence, in view of
\eqref{rem32} then ${M}^{\varrho}\simeq N$, and thus 
${M}^{\varrho}$  especially obeys ${M}^{\varrho}\not={\mathbb{C}}\,{\mathbf 1}$.
\item\label{rem34}
Note in context of \eqref{add.2} that if $\varOmega$ is separating for a
$vN$-algebra $N$, then according to \cite{Arak:72}
the condition $\varphi\in \overline{N_+\varOmega}$ is equivalent to 
the existence of a densely defined, positive, selfadjoint linear operator $A$
which is {\em affiliated} with
$N$, $A\,{\mathop{\eta}}\, N$, and which obeys $\varphi=A\varOmega$. On the other hand
one knows that the latter is also equivalent to  $\langle (\cdot)|_{N^{\,\prime}}\varphi,\varOmega\rangle\geq 0$.
\item\label{rem35}
Especially, if for faithful $\varrho$ a cyclic and separating
$\varOmega\in {\mathcal S}_M(\varrho)$ is chosen, then owing to
$\langle(\cdot)|_{M^{\,\prime}}\psi_\varOmega^\varrho(\nu),\varOmega\rangle\geq 0$,
the previous facts ensure existence of a unique positive selfadjoint linear operator
$T \,{\mathop{\eta}}\,M$ with $\psi_\varOmega^\varrho(\nu)=T\varOmega$.
In view of the well-known Radon-Nikodym theorem \cite{Saka:66} and following \cite{Arak:72}
$T$ then is referred to as generalized Sakai's Radon-Nikodym operator of $\nu$
with respect to $\varrho$. In this case $j(\nu|\varrho)$ of Remark \ref{rem2}\,\eqref{rem2.3} then
may be written as
$
j(\nu|\varrho)=\frac{1}{2}\,
\|J_\varOmega T^* J_\varOmega\varOmega-T\varOmega\|^2\,.
$
Thus, in the case at hand $j(\nu|\varrho)$ amounts to be the same as the {\em skewinformation}
$I(\varrho,T)$ of the generalized Radon-Nikodym
operator $T$ with respect to $\varrho$ in the sense of \cite{CoSt:78}, see
\cite{ArYa:60,WiYa:63,WiYa:64,Yana:64} for background informations, and which together with
the other facts from Remark \ref{rem2}\,\eqref{rem2.3} might facilitate the understanding
of the notion of commutation also in those cases which go beyond Example \ref{ex3}, see also 
\cite[Definition 3,\,Lemma 2,\,Theorem 2]{Albe:97.1}.
\end{enumerate}
\end{remark}
Close this part by two auxiliary results, which are versions of more 
general assertions, and which  
ensure applicability of Lemma \ref{cent}/\ref{centfac}, in our cases of interest. 
In the following, for each   
$a\in M$ and $\mu\in M_{*+}$, we let a  
normal positive linear form $\mu^a$ be defined by $\mu^a(x)=\mu(a^*x a)$, for all $x\in M$.
\begin{lemma}\label{premain0}
Let $M$ be a standard form $vN$-algebra, and be $\varrho\in M_{*+}$.
Suppose $\nu=\varrho^x$, for some $x\in M_+$. Then $\psi_\varOmega^\varrho(\nu)=x\xi_\varrho$ is fulfilled. 
Moreover, in case that $\varrho$ is faithful, $\nu$ is commuting with $\varrho$ 
if, and only if, $x\in \bigl(M^\varrho\bigr)_+$ holds.
\end{lemma}
\begin{proof}
Let $\varOmega$ be a cyclic and separating vector for $M$. 
Note that by assumption on $\nu$, $x\xi_\varrho \in {\mathcal S}_M(\nu)$. Also, by positivity of 
$x$, $h=h_{x\xi_\varrho,\xi_\varrho}=h_{\sqrt{x}\xi_\varrho,\sqrt{x}\xi_\varrho}\geq 0$. 
From this then $s(h)=p^{\,\prime}(\sqrt{x}\xi_\varrho)\leq p^{\,\prime}(x\xi_\varrho)$ 
follows. On the 
other hand, in any case, $p^{\,\prime}(x\xi_\varrho)\leq p^{\,\prime}(\sqrt{x}\xi_\varrho)$ has 
to be fulfilled. Thus we conclude  
$s(h)=p^{\,\prime}(x\xi_\varrho)$. In view of 
Theorem \ref{bas5}\,\eqref{bas53} and 
Lemma \ref{sesq}\,\eqref{positiv0*} then $d_{\mathrm{B}}(\nu,\varrho)=\|x\xi_\varrho-\xi_\varrho\|$ and 
$\nu^\perp=0$ follow. In the case at hand Lemma \ref{mod1}\,\eqref{red.3} can be applied  
and yields $\psi_\varOmega^\varrho(\nu)=x\xi_\varrho$. 

Suppose $\varrho$ is faithful. Then we can choose a cyclic and separating vector
$\varOmega\in {\mathcal S}_M(\varrho)$, and will work in the corresponding standard form action of $M$. 
Suppose that $\nu$ commutes with $\varrho$. 
By Definition \ref{commu}, since $\xi_\varrho=\varOmega$ holds, and since 
commutation means $\psi_\varOmega^\varrho(\nu)=\xi_\nu$, we have    
$\langle (\cdot)|_{M^{\,\prime}}\,\xi_\nu,\varOmega\rangle=h_{\xi_\nu,\xi_\varrho}\geq 0$.
But in view of Lemma \ref{rem2.1a}\,\eqref{rem22} and
Lemma \ref{zentral} also $\langle (\cdot)|_{M}\xi_\nu,\varOmega\rangle\geq 0$.
Hence, $\xi_\nu\in \overline{M_+\varOmega}\cap \overline{M_+^{\,\prime}\varOmega}$ by
Remark \ref{rem3}\,\eqref{rem34}. As has been explained in context of \eqref{add.2}, 
this is the same as $\xi_\nu\in \overline{(M^\varrho)_+\varOmega}$.
Also, since
$M^\varrho\subset M$ implies $M^{\,\prime}\subset \bigl(M^\varrho\bigr)^{\,\prime}$,
$\varOmega$ must be separating for $M^\varrho$.
Hence, by Remark \ref{rem3}\,\eqref{rem34} and with $N=M^\varrho$, $\xi_\nu=A\varOmega$ follows,
with densely defined,
positive, selfadjoint linear operator $A$ which is affiliated with $M^\varrho\subset M$. 
In view of the above then $\psi_\varOmega^\varrho(\nu)=x\varOmega=A\varOmega$. 
Since $\varOmega$ is separating, from 
this $x=A$ follows, that is, $x\in \bigl(M^\varrho\bigr)_+$. 

On the other hand, for $x\in M^\varrho$, $x\geq {\mathbf 0}$, in view of the above 
formula $\psi_\varOmega^\varrho(\nu)=x\varOmega$ and \eqref{add.1}, $\psi_\varOmega^\varrho(\nu)
\in {\mathcal H}_{\mathrm{sa}}^\varOmega$ follows. Applying Lemma \ref{mod1}\,\eqref{red.2} 
accordingly yields $\psi_\varOmega^\varrho(\nu)\in {\mathcal P}_\varOmega^\natural$. 
Hence, $\nu$ commutes with $\varrho$. 
\end{proof}
\begin{remark}\label{premain1}
According to Remark \ref{rem3}\,\eqref{rem34}, if $\varrho$ is faithful, 
and $\varOmega\in {\mathcal S}_M(\varrho)$ 
is a cyclic and separating vector, then for each $\nu\in M_{*+}$ one has 
$\psi_\varOmega^\varrho(\nu)=x\varOmega$, with a densely defined, affiliated 
with $M$, selfadjoint positive linear operator $x$. Thus, by literally 
the same conclusions as in the previous proof we see that    
$\nu$ will be commuting with $\varrho$ iff $x=x^*\geq {\mathbf 0}$ is affiliated 
with $M^\varrho$. 
\end{remark}
\begin{lemma}\label{premain}
Let $M$ be a standard form $vN$-algebra, and be $\varrho\in M_{*+}$ faithful.
For each $\nu\in M_{*+}$ commuting with $\varrho$ one has $s(\nu)\in M^\varrho$.
\end{lemma}
\begin{proof}
Let us consider the standard form action of $M$ with respect to cyclic and separating
$\varOmega\in {\mathcal S}_M(\varrho)$. By assumption such a vector has to exist. 
As mentioned in the previous proof, $\psi_\varOmega^\varrho(\nu)=\xi_\nu$ implies  
(in fact is equivalent to) $\xi_\nu=A\varOmega$, 
with densely defined,
positive, selfadjoint linear operator $A$ which is affiliated with $M^\varrho$.
Thus $s(\nu)\leq q$, where $q$ is the orthoprojection
onto the closure of the range of $A$. By spectral calculus, if
$\{E(\lambda)\}\subset M^\varrho$ is the (left-continuous) spectral resolution of $A$,
then $q={\mathrm{l.u.b.}}\{E(\lambda)-E(0+):\lambda>0\}\in M^\varrho$. Also, owing to
$AE(\lambda)\in (M^\varrho)_+$, for $\nu_\lambda$
with $E(\lambda)\xi_\nu=AE(\lambda)\varOmega\in {\mathcal S}_M(\nu_\lambda)$,
we have $\nu_\lambda\leq \nu$ and $s(\nu_\lambda)=E(\lambda)-E(0+)$. Hence 
$q\leq s(\nu)$. Thus $q=s(\nu)$ follows. In view of the above then $s(\nu)\in M^\varrho$.
\end{proof}
\subsubsection{Conclusions in the standard form case}\label{finis}
Start with a mechanism for creating pairs of positive linear forms 
which are relevant in context of Problem \ref{prob}\,\eqref{prob4}, and   
which have been found useful in other context \cite[{\sc Lemma 4,Theorem 6}]{AlUh:00.1}, too.
\begin{lemma}\label{main0}
Let $M$ be an infinite standard form $vN$-algebra, and be $\varrho\in M_{*+}$ faithful. 
For each  $\sigma\in M_{*+}$ with $s(\sigma)\sim {\mathbf 1}$ there is $a\in M$ 
with ${\mathcal S}_M(\sigma^a|\varrho)\not={\mathcal S}_M(\sigma^a)$. In particular, 
there is an orthoprojection $z\in M^\varrho$ such that $M\simeq zMz$ holds, and 
${\mathcal S}_M(\mu|\varrho)\not={\mathcal S}_M(\mu)$ is fulfilled, for each 
$\mu\in M_{*+}$ 
with support $s(\mu)=z$.
\end{lemma}
\begin{proof}
Let $M=Mc+Mc^\perp$, with
orthoprojection $c\in M\cap M^{\,\prime}$, be the canonical
decomposition of $M$ into a finite component $Mc$ (which might be vanishing) and properly
infinite component $Mc^\perp$. By central decomposition techniques, the properly infinite 
case can be reduced to an appropriately defined field $\{M_\lambda\}$ of properly infinite `subfactors'
$M_\lambda$ acting over direct integral Hilbert subspaces ${\mathcal H}_\lambda\subset c^\perp{\mathcal H}$, in the sense that each $x\in Mc^\perp$ and $\omega\in (Mc^\perp)_{*+}$ can be written as appropriately defined central integrals $x=\int
\oplus x_\lambda d\/\mu(\lambda)$ and $\omega=\int
\oplus\, \omega_\lambda d\/\mu(\lambda)$, with some $x_\lambda\in M_\lambda$ and 
$\omega_\lambda\in (M_\lambda)_{*+}$, such that $\omega(x)$ can be written as
$\omega(x)=\int
 \omega_\lambda(x_\lambda) d\/\mu(\lambda)$  over terms $\omega_\lambda(x_\lambda)$, with 
$\omega_\lambda=\omega|_{M_\lambda}$, and  
with respect to some central measure $\mu(\lambda)$. Thereby, the index $\lambda$ refers to the
central decomposition which is isometrically defined relative to some direct integral decomposition 
$c^\perp{\mathcal H}=\int\oplus{\mathcal H}_\lambda d\/\mu(\lambda)$ of the Hilbert space 
$c^\perp{\mathcal H}$,   
see e.g.~\cite[{\em Chapter I.},\,5.,\,10.~Corollary]{Schw:67} for the details. 
Since all decompositions refer to the center and $\varrho$ is faithful over $M$, 
for $\varrho_\lambda=\varrho|_{M_\lambda}$ we will 
have $(M^\varrho c^\perp)_\lambda=M_\lambda^{\varrho_\lambda}\subset M_\lambda$, with 
faithful $\varrho_\lambda$ over properly infinite subfactor $M_\lambda$. 
By Remark \ref{rem3}\,\eqref{rem33}, $M_\lambda^{\varrho_\lambda}$ cannot be trivial. 
Therefore, there is an orthoprojection ${\mathbf 0}<c_\lambda<{\mathbf 1}_\lambda$, with 
$c_\lambda\in (M^\varrho c^\perp)_\lambda$. As orthoprojection of the infinite factor  
$M_\lambda$, at least one of $c_\lambda$ and $c_\lambda^\perp$ must be infinite. We make a choice, and 
call this infinite projection 
$z_\lambda$. Since $M$ is countably decomposable by assumption, each of $M_\lambda$ is a 
countably decomposable infinite factor. Hence,  
we then must have 
$z_\lambda<{\mathbf 1}_\lambda$ and $z_\lambda\sim {\mathbf 1}_\lambda$, with respect to  $M_\lambda$. 
Define an orthoprojection $z$ in $M$ as follows:
$$
z=c+\int \oplus z_\lambda d\/\mu(\lambda)\,.
$$
One then has $z\in M^\varrho$, $z<{\mathbf 1}$, and
$z\sim {\mathbf 1}$ with respect to $M$. But then, by assumption on $\sigma$, 
there is $a\in M$ with $a^*a=s(\sigma)$ and $aa^*=z$. 
It is easily seen that for each such $a$ then $s(\sigma^a)=z$ must be fulfilled. 
Hence, $s(\sigma^a)\in M^\varrho$ with 
$s(\sigma^a)<{\mathbf 1}=s(\varrho)$ but $s(\sigma^a)\sim {\mathbf 1}$, for each such $a$. 
By Lemma \ref{cent} with $\nu=\sigma^a$ then 
${\mathcal S}_M(\sigma^a|\varrho)\not={\mathcal S}_M(\sigma^a)$ follows, for each $a\in M$ with $a^*a=s(\sigma)$ and $aa^*=z$. 

Now, let $v\in M$ be a partial isometry with $v^*v=z$ and $vv^*={\mathbf 1}$. Then, the map 
$\pi:M\ni x\mapsto v^*xv\in zMz$ is a $^*$-isomorphism between $M$ and the hereditary $vN$-subalgebra $zMz$. Thus, for 
$\mu$ running through all $\mu\in M_{*+}$ with $s(\mu)=z$, $\sigma=\mu^v$ is running 
through all faithful $\sigma\in M_{*+}$. In view of $s(\sigma)={\mathbf 1}$, by the above 
conclusion for each such $\sigma$ and with $a=v^*$ the result is  
${\mathcal S}_M(\sigma^a|\varrho)\not={\mathcal S}_M(\sigma^a)$. 
But note that $\sigma^a=\sigma(v(\cdot)v^*)=\mu^v(v(\cdot)v^*)=\mu(v^*v(\cdot)v^*v)=\mu(z(\cdot)z)=
\mu$. That is, ${\mathcal S}_M(\mu|\varrho)\not={\mathcal S}_M(\mu)$ must be fulfilled, for each 
$\mu\in M_{*+}$ with $s(\mu)=z$.
\end{proof}
\begin{remark}\label{main2}
Note that under the premises of Lemma \ref{main0} there have to exist both,  
$\nu\in M_{*+}$ which commute/do not commute with $\varrho$, but which obey 
${\mathcal S}_M(\nu|\varrho)\not={\mathcal S}_M(\nu)$. In fact, according to 
Lemma \ref{premain0} and  Lemma \ref{main0}, $\nu=\varrho^z$ provides a special 
positive linear form  
commuting with $\varrho$ and obeying  
${\mathcal S}_M(\nu|\varrho)\not={\mathcal S}_M(\nu)$. But according to Lemma \ref{main0} 
the latter remains true also 
for any other $\nu=\varrho^x$, with invertible (within $zMz$) 
$x\in (zMz)_+$. If each of these functionals were commuting with $\varrho$ then, by  
Lemma \ref{premain0}, $x\in M^\varrho$  
for all those $x$. By taking the closure we had $(zMz)_+\subset M^\varrho$, 
and therefore $zMz\subset  M^\varrho$. 
But then $\varrho|_{zMz}$ would be a faithful tracial 
normal positive linear form over $zMz$. This is a contradiction, since owing to  
$zMz\simeq M$, $zMz$ cannot be finite. 
\end{remark}
\begin{example}\label{finrem}
If $M$ is a factor of type $\mathrm{I}_\infty$, things around Lemma \ref{main0} 
read quite explicit. In fact, in terms
of Example \ref{ex2}/\ref{ex3} we then have $\varrho=\tau_c$, with
uniquely determined, positive definite 
$\tau$-trace class operator $c$. By Example \ref{ex3}, in this case 
$M^\varrho=\{x\in M: xc=cx\}$. 
Then, orthoprojections $z\in M$ obeying $zc=cz$, 
$\tau(z)=\infty$ and $z\not={\mathbf 1}$ must exist and, by   
Lemma \ref{centfac} and in line with the last part of the proof of Lemma \ref{main0}, this is the set of all 
orthoprojections $z$ with respect 
to which the hypothesis of Lemma \ref{main0} holds,  
and therefore each $\nu\in M_{*+}$ with support in this set provides an example with    
${\mathcal S}_M(\nu|\varrho)\not={\mathcal S}_M(\nu)$. 
\end{example}
In view of Problems \ref{prob}\,\eqref{prob3}/\eqref{prob4}, we may summarize all that also 
as follows.
\begin{theorem}\label{main}
Let $\{M,\varOmega\}$ be a standard form $vN$-algebra, $\varrho\in M_{*+}$ faithful.
\begin{enumerate}
\item\label{main.1}
${\mathcal S}_M(\nu|\varrho)={\mathcal S}_M(\nu)$ holds
for all $\nu\in M_{*+}$ if, and only if, $M$ is finite.
\item\label{main.2}
For infinite  $M$, there exists $\nu\in M_{*+}$ obeying 
\begin{equation}\label{commutation}
d_{\mathrm{B}}(\nu,\varrho)=\|\xi_\nu-\xi_\varrho\|
\end{equation}
such that ${\mathcal S}_M(\nu|\varrho)\not={\mathcal S}_M(\nu)$.
\item\label{main.3}
For $\nu\in M_{*+}$ satisfying \eqref{commutation}, ${\mathcal S}_M(\nu|\varrho)=
\{u\xi_\nu:\,u^*u={\mathbf 1},\,u\in M^{\,\prime}\}$ holds.
\end{enumerate}
\end{theorem}
\begin{proof}
By Lemma \ref{mincommu}, condition \eqref{commutation} means commutation of $\nu$ with $\varrho$, 
in the sense of Definition \ref{commu}. Thus \eqref{main.3} follows along with Example 
\ref{main1} and formula \eqref{allg}. Relating \eqref{main.1}, by Remark \ref{fin} we yet
know that ${\mathcal S}_M(\nu|\varrho)={\mathcal S}_M(\nu)$ is
always fulfilled on a finite $M$. On the other hand, for infinite $M$, according to  
Lemma \ref{main0} and Remark \ref{main2} there exists $\nu$ commuting with $\varrho$ and obeying  
${\mathcal S}_M(\nu|\varrho)\not={\mathcal S}_M(\nu)$, which in view of Lemma \ref{mincommu} 
yields \eqref{main.2}.  
Since a $vN$-algebra is
either finite or infinite, this then also completes the proof of \eqref{main.1}. 
\end{proof}
In the factor case, for mutually commuting normal positive linear forms one of which is faithful,
things considerably simplify (see Example \ref{finrem}), and then 
Theorem \ref{main}\,\eqref{main.2} can be supplemented by 
the following detail.
\begin{lemma}\label{mainend}
Let $M$ be a factor acting in standard form. 
Suppose $\nu\in M_{*+}$ obeys condition \eqref{commutation} with faithful $\varrho\in M_{*+}$.
Then, ${\mathcal S}_M(\nu|\varrho)={\mathcal S}_M(\nu)$ holds if, and only if,
$\nu$ is faithful or the support orthoprojection $s(\nu)$ is finite.
\end{lemma}
\begin{proof}
By Lemma \ref{mincommu} and Lemma \ref{premain}, condition \eqref{commutation} implies 
$s(\nu)\in M^\varrho$. Lemma \ref{centfac} then may be applied, and yields 
${\mathcal S}_M(\nu|\varrho)={\mathcal S}_M(\nu)$ iff either $s(\nu)$ is finite or
$s(\nu)={\mathbf 1}$, cf.~also Remark \ref{rem42}.
\end{proof}
\begin{example}\label{finrem1}
For a type-\,$\mathrm{III}$\,-factor $M$, things around commutation 
are less explicit than in the type\,-$\mathrm{I}_\infty$-\,case. A characterization like \eqref{commutation} 
which reads in terms of a geometrical condition in context of the Bures distance then possibly   
is the best one can do. By  
Lemma \ref{mincommu}, Theorem \ref{main}\,\eqref{main.2} and 
Lemma \ref{mainend} tell us that for faithful $\varrho$ there is      
$\nu\in M_{*+}$ obeying \eqref{commutation} and    
${\mathcal S}_M(\nu|\varrho)\not={\mathcal S}_M(\nu)$, and the latter 
remains true for any $\nu$ with \eqref{commutation}, unless either $\nu$ is faithful or is vanishing.
\end{example}

\section*{Acknowledgments}
We gratefully remember discussions with former students U.~Meister,
S.~Schwebel and U.~Steinmetz, and thank for
contributing to the `Project Seminar on Bures Geometry' at Mathematical Institute,
in an earlier stage of the project.
Thanks go to `Deutsche Forschungsgemeinschaft'
for funding the research of one of us (PMA) by grants. In particular, 
thanks go to Eberhard Zeidler (MPI MIS Leipzig), Gerd Rudolph
(ITP Leipzig) and the Institute of Theoretical Physics (ITP) for supporting the
research program, and for kind hospitality. Last but not least, we have to
thank Bodo Geyer, speaker of the `Graduiertenkolleg Quantenfeldtheorie', for supporting both,
lecturing (PMA) and graduating (GP) at the Graduiertenkolleg.

\bibliographystyle{abbrv}

\end{document}